\documentclass[12pt,a4paper]{article}
\usepackage{amsmath,amssymb}
\usepackage[english]{babel}
\usepackage{eucal}

\DeclareMathVersion{scriptstyle}

\usepackage{isolatin1,theorem}
\theoremstyle{change}
{\theorembodyfont{\slshape} 
  \newtheorem{thm}{Theorem.}[section]
  \newtheorem{lemma}[thm]{Lemma.}
  \newtheorem{prop}[thm]{Proposition.}
  \newtheorem{cor}[thm]{Corollary.}
}
{\theorembodyfont{\rmfamily} 
  \newtheorem{mastereq}[thm]{Theorem. (Master equation)}
  \newtheorem{masterineq}[thm]{Theorem. (Master inequality)}
  
  \newtheorem{remark}[thm]{Remark.}

}

\numberwithin{equation}{section} 

\newcommand{\proof}[1][Proof. ]{{\it#1}}

\def\endproof{{\nobreak\qquad{\scriptstyle \blacksquare}}}

\def\C{{\mathbb C}}
\def\B{{\mathbb B}}
\def\E{{\mathbb E}}

\def\H{{\mathbb H}}
\def\N{{\mathbb N}}
\def\R{{\mathbb R}}
\def\Y{{\mathbb Y}}

\def\V{{\mathbb V}}

\def\CA{{\mathcal A}}
\def\CB{{\mathcal B}}
\def\CC{{\mathcal C}}
\def\CD{{\mathcal D}}

\def\CF{{\mathcal F}}

\def\CJ{{\mathcal J}}
\def\CK{{\mathcal K}}
\def\CL{{\mathcal L}}
\def\CM{{\mathcal M}}

\def\CO{{\mathcal O}}

\def\CU{{\mathcal U}}

\def\unit{{\bf 1}}
\def\i{{\rm i}}
\def\e{{\rm e}}
\def\d{{\rm d}}
\def\eps{\varepsilon}
\def\<{{\langle}}
\def\>{{\rangle}}
\def\lmax{{\lambda_{\max}}}
\def\lmin{{\lambda_{\min}}}
\def\id{{\text{\rm id}}}
\def\tr{{\text{\rm tr}}}
\def\Tr{{\text{\rm Tr}}}
\def\GRM{{\text{\rm GRM}}}
\def\SGRM{{\text{\rm SGRM}}}
\def\GRMR{{\text{\rm GRM}}^{\R}}
\def\GOE{{\text{\rm GOE}}}
\def\GOES{{\text{\rm GOE}}^\ast}
\def\GSE{{\text{\rm GSE}}}
\def\GSES{{\text{\rm GSE}}^\ast}
\def\im{{\text{\rm Im}}}
\def\re{{\text{\rm Re}}}
\def\supp{{\text{\rm supp}}}
\def\grad{{\text{\rm grad}}}
\def\diff{\frac{\rm d}{{\rm d}t}\Big|_{t=0}}
\def\diffz{\frac{\rm d}{{\rm d}z}\Big|_{z=0}}
\def\diffl{\frac{\rm d}{{\rm d}\lambda}}

\def\cc{^*}
\def\Om{\CO^{(-)}}

\def\Distc{\CD_c'(\R)}
\def\Cinf{C^\infty(\R)}
\def\Ccinf{C_c^\infty(\R)}

\def\unitH{\begin{pmatrix}
1 & 0\\
0 & 1
\end{pmatrix}}
\def\j{\begin{pmatrix}
\i & 0\\
0 & -\i
\end{pmatrix}}
\def\k{\begin{pmatrix}
0 & 1\\
-1 & 0
\end{pmatrix}}
\def\l{\begin{pmatrix}
0 & \i\\
\i & 0
\end{pmatrix}}

\begin{document}

\title{Non-commutative Polynomials of \\ Independent Gaussian Random Matrices.\\ The Real and Symplectic Cases.}
\author{\textsc{Hanne Schultz}\footnote{This work was partially supported
    by MaPhySto -- A Network in Mathematical Physics and Stochastics,
    funded by The Danish National Research Foundation.} \footnote{ As a
    student of the PhD-school OP-ALG-TOP-GEO the author is partially
    supported by the Danish Research Training Council.}}
\date{}       
\maketitle

\begin{abstract} \noindent In [HT2] Haagerup and Thorbj\o rnsen prove the
  following extension of Voiculescu's random matrix model (cf. [V2,
  Theorem~2.2]): For each $n\in\N$, let $X_1^{(n)}, \ldots, X_r^{(n)}$ be
  a system of $r$ independent complex self-adjoint random matrices from the
  class $\SGRM(n, \frac 1n)$, and let $x_1, \ldots, x_r$ be a semicircular
  system in a $C\cc$-probability space. Then for any
  polynomial $p$ in $r$ non-commuting variables the convergence
  \[
  \lim_{n\rightarrow \infty}\|p(X_1^{(n)}, \ldots, X_r^{(n)})\| = \|p(x_1,
  \ldots, x_r)\|
  \]
  holds almost surely. We generalize this result to sets of independent
  Gaussian random matrices with real or symplectic entries (the $\GOE$- and
  the $\GSE$-ensembles) and random matrix ensembles related to these.
\end{abstract}
  
\section{Introduction.}
 
Throughout this paper we let $(\Omega, \CF, P)$ denote a fixed probability
space, and for each positive integer $n$ and each $\sigma > 0$ we
let $\SGRM(n, \sigma^2)$ denote the set of $n\times n$ self-adjoint
Gaussian random matrices defined in \cite{HT2}:
\begin{itemize}
  \item[(i)] $\SGRM(n, \sigma^2)$ is the set of self-adjoint random
  matrices $X = (X_{ij}): \Omega \rightarrow M_n(\C)$ satisfying that
  $\{X_{ii}|1\leq i\leq n\}\cup \{\sqrt 2 \re X_{ij}|1\leq i<j \leq n\}
  \cup \{\sqrt 2 \im X_{ij}|1\leq i<j \leq n\}$ is a set of $n^2$
  i.i.d. random variables with distribution $N(0,  \sigma^2)$.
\end{itemize}

We shall also consider the following related random matrix ensembles:

\begin{itemize}
  \item[(ii)] $\GRM(n, \sigma^2)$, which was defined in [HT2] as well, is the set of random matrices $Y : \Omega
  \rightarrow M_n(\C)$ satisfying that the real and the imaginary parts of
  the entries of $Y$, $\re(Y_{ij}),\; \im (Y_{ij}),\; 1\leq i, j\leq n$,
  constitute a set of 2$n^2$ i.i.d. random variables with distribution $N(0, \frac 12  \sigma^2)$.
  \item[(iii)]  $\GRMR(n,  \sigma^2)$ is the set of random matrices $Y :
  \Omega \rightarrow M_n(\mathbb{R})$ satisfying that the entries of $Y$, $Y_{ij}\;1\leq i, j\leq n$, constitute a set of $n^2$ i.i.d. random variables with distribution $N(0,  \sigma^2)$.
  \item[(iv)] $\GOE(n,  \sigma^2)$ is the set of self-adjoint real random
  matrices $X : \Omega \rightarrow M_n(\R)$ satisfying that
  $\{\frac{1}{\sqrt 2}X_{ii}|1\leq i \leq n\}\cup\{X_{ij}|1\leq i <j \leq
  n\} $ is a set of $\frac 12 n(n+1)$ i.i.d. random variables with distribution $N(0,  \sigma^2)$.
  \item[(v)] $\GOES(n, \sigma^2)$ is the set of self-adjoint purely
  imaginary random matrices  $X : \Omega \rightarrow M_n(\i\R)$ satisfying
  that $\{\im (X_{ij})|1\leq i <j \leq n\} $ is a set of $\frac 12 n(n-1)$ i.i.d. random variables with distribution $N(0,  \sigma^2)$.  
\end{itemize}

Note that (up to scaling) $\SGRM(n, \sigma^2)$ is
the Gaussian Unitary Ensemble (${\rm GUE}$) from [Me], and $\GOE(n,
\sigma^2)$ is the Gaussian Orthogonal Ensemble ($\GOE$) from [Me].

Also note that for every $Y\in \GRMR(n,  \sigma^2)$ we have 
\[
X_1=\frac{1}{\sqrt{2}} (Y + Y^*)  \in  \GOE(n,  \sigma^2),
\] 
and 
\[
X_2=\frac{\textrm{1}}{\i\sqrt{2}} (Y - Y^*) \in \GOES(n,  \sigma^2).
\]
Moreover, $X_1$ and $X_2$ are stochastically independent. Conversely, if $X_1\in \GOE(n,  \sigma^2)$ and $X_2\in  \GOES(n,  \sigma^2)$ are independent, then 
\begin{equation}\label{GRMR}
Y = \frac{1}{\sqrt 2}(X_1+\i X_2)\in \GRMR(n, \sigma^2).
\end{equation}

In their paper, \emph{A new application of random matrices: ${\rm Ext}(C_{red}\cc (F_2))$ is not a group}, Haagerup and Thorbj\o rnsen prove

{\bf Theorem.} [HT2, Theorem~7.1]
Let $r \in \N$, and for each $n\in \N$, let $X_1^{(n)}, \ldots, X_r^{(n)}$ be stochastically independent random matrices from $\SGRM(n,\frac 1n)$. Furthermore, let $(\CA, \tau)$ be a $C\cc$-probability space with $\tau$ a faithful state on $\CA$, and let $\{x_1, \ldots, x_r\}$ be a semicircular system in $(\CA, \tau)$. Then there is a $P$-null set $N\subseteq\Omega$ such that for any $\omega\in \Omega\setminus N$ and for any polynomial $p$ in $r$ non-commuting variables:
\[
\lim_{n\rightarrow \infty}\|p(X_1^{(n)}(\omega), \ldots, X_r^{(n)}(\omega))\| = \|p(x_1, \ldots, x_r)\|.
\]

\vspace{.2cm}

We prove similar results for some other classes of self-adjoint Gaussian random matrices:

{\bf Theorem~A.}
Let $r$, $s \in \N_0$ with $r+s\geq 1$ , and for each $n\in \N$, let $X_1^{(n)}, \ldots, X_{r+s}^{(n)}$ be stochastically independent random matrices defined on $(\Omega, \CF, P)$ such that $ X_1^{(n)}, \ldots,  X_r^{(n)}\in \GOE(n, \frac 1n)$ and $ X_{r+1}^{(n)},\ldots, X_{r+s}^{(n)} \in \GOES(n, \frac 1n)$. Furthermore, let $(\CA, \tau)$ be a $C\cc$-probability space with $\tau$ a faithful state on $\CA$, and let $\{x_1, \ldots, x_{r+s}\}$ be a semicircular system in $(\CA, \tau)$. Then there is a $P$-null set $N\subseteq\Omega$ such that for any $\omega\in \Omega\setminus N$ and for any polynomial $p$ in $r+s$ non-commuting variables:
\[
\lim_{n\rightarrow \infty}\|p(X_1^{(n)}(\omega), \ldots, X_{r+s}^{(n)}(\omega))\| = \|p(x_1, \ldots, x_{r+s})\|.
\]

\vspace{.2cm}

In the sections 2 to 5 we concentrate on proving that there is a $P$-null set $N' \subseteq \Omega$ such that for any non-commutative polynomial $p$ in $r+s$ variables and for every $\omega \in \Omega \setminus N'$,
\begin{equation}\label{B}
  \limsup_{n \rightarrow \infty}\|p(X_1^{(n)}(\omega), \ldots , X_{r+s}^{(n)}(\omega))\| \leq \|p(x_1, \ldots , x_{r+s})\|,
\end{equation}

and in Section 6 we prove that there is a $P$-null set $N'' \subseteq \Omega$ such that for any non-commutative polynomial $p$ in $r+s$ variables and for every $\omega \in \Omega \setminus N''$,
\begin{equation}\label{C}
  \liminf_{n \rightarrow \infty}\|p(X_1^{(n)}(\omega), \ldots , X_{r+s}^{(n)}(\omega))\| \geq \|p(x_1, \ldots , x_{r+s})\|.
\end{equation}

\vspace{.2cm}

Clearly, Theorem~A follows from (\ref{B}) and (\ref{C}). 

The proof of (\ref{B}) follows the lines of \cite{HT2}, and we shall at some places leave out the details and refer to that paper instead. However, additional difficulties arise in the $\GOE / \GOES$-case, and this is
mainly due to the appearance of a term of order $\frac 1n$ in our
''master equation'' (cf. Theorem~\ref{MEQ}),
\[
 \E\Big\{(a_0-\lambda)H_n(\lambda) + \sum_{j=1}^r[a_j H_n(\lambda)a_j^*
 H_n(\lambda) + a_j^*H_n(\lambda)a_j H_n(\lambda)] + \mathbf{1}_m\Big\} =
 -\frac{1}{n}R_n(\lambda),
 \]
where we use the notation introduced in Section~2. The corresponding equation in \cite{HT2} (cf. [HT2, Theorem~3.6]) contains
no such term:
\[
\E\Big\{(a_0-\lambda)H_n(\lambda) + \sum_{j=1}^r a_jH_n(\lambda)a_jH_n(\lambda)
+ \unit_m\Big\} = 0.
\]
A crucial ingredient in the proof given by Haagerup and Thorbj\o rnsen is the following estimate (cf. [HT2, Theorem~5.7]):
\[
\|G_n(\lambda)-G(\lambda)\| = O\Big(\frac{1}{n^2}\Big),
\]
Because of the difference mentioned above we get an extra term of order $\frac 1n$ (cf. Theorem~\ref{longthm}):
\[
\Big\|G_n(\lambda)-G(\lambda)-\frac 1n L(\lambda)\Big\| =  O\Big(\frac{1}{n^2}\Big).
\]

In Section~5 it is proved that $l(\lambda) := \tr_m(L(\lambda\unit_m))$
gives rise to a compactly supported distribution, $\Lambda$, such that 
\[
\E\{(\tr_m\otimes\tr_n)\phi(S_n)\}=(\tr_m\otimes\tau)\phi(s)+\frac 1n
\Lambda(\phi) +  O\Big(\frac{1}{n^2}\Big)
\]
for any $\phi \in \Ccinf$ (still with the same notation as in Section~2).  
In fact, $\supp(\Lambda)\subseteq \sigma(s)$, so $\Lambda(\phi) = 0$ for any $\phi \in \Ccinf$ with $\supp(\phi)\cap \sigma(s) = \emptyset$. It follows that for  $\psi \in \Cinf$ such that $\psi$ is constant outside a compact set of $\R$ and $\supp(\psi)\cap \sigma(s) = \emptyset$,
\[
\E\{(\tr_m\otimes\tr_n)\psi(S_n)\} =  O\Big(\frac{1}{n^2}\Big)
\]
- an estimate similar to the one obtained by Haagerup and Thorbj\o rnsen in [HT2, Lemma~6.3] and a cornerstone of the proof of Theorem~A.

\vspace{.2cm}

In Section~7 we consider yet two other ensembles which are random matrix
ensembles with quaternionic entries. Remember that the quaternions is the division ring, $\H$, which is, as a vector space over the real numbers, spanned by four linearly independent vectors, $\unit, j, k, l$, satisfying the identities
\[
j^2 = k^2 = l^2 = -\unit 
\]
and
\[
jk = -kj = l, \qquad kl = -lk = j, \qquad lj = -jl = k.
\]
We denote by $\H^\C$ the complexification of $\H$.

It is well known that 
\[
\unit \mapsto \unitH, \qquad j \mapsto \j, \qquad k \mapsto \k, \qquad l
\mapsto \l
\]
defines a ring homomorphism which is an embedding of $\H$ into $M_2(\C)$. 
\begin{itemize}
\item[(vi)] $\GSE(n,\sigma^2)$ is the set of random matrices $X\rightarrow
  M_n(\H)= \H\otimes_\R M_n(\R)$ satisfying that
\[
X = \unit\otimes V + j\otimes (\i W) + k\otimes (\i Y) + l\otimes (\i Z),
\]
for some $V\in \GOE(n, \frac{\sigma^2}{4})$ and some $W, Y, Z\in
\GOES(n,\frac{\sigma^2}{4})$, where $V, W, Y$ and $Z$ are stochastically
independent. If we identify $\H$ with a real sub-algebra of $M_2(\C)$ in the
way described above, then
\[
X = \unitH\otimes V + \j\otimes (\i W) + \k\otimes (\i Y) + \l\otimes (\i Z),
\]
as a random matrix taking values in $M_2(\C)\otimes M_n(\C)$, and this
shall be our preferred description of $\GSE(n,\sigma^2)$.
\item[(vii)] $\GSES(n,\sigma^2)$ is the set of random matrices $X\rightarrow
  M_n(\H^\C)= \H\otimes_\R M_n(\C)$ satisfying that
\[
X = \unit\otimes V + j\otimes (\i W) + k\otimes (\i Y) + l\otimes (\i Z),
\]
for some $V\in \GOES(n, \frac{\sigma^2}{4})$ and some $W, Y, Z\in
\GOE(n,\frac{\sigma^2}{4})$, where $V, W, Y$ and $Z$ are stochastically
independent. Again, we shall preferably consider $X$ as a
random matrix taking values in $M_2(\C)\otimes M_n(\C)$ with
\[
X = \unitH\otimes V + \j\otimes (\i W) + \k\otimes (\i Y) + \l\otimes (\i Z).
\]
\end{itemize}

Up to scaling $\GSE(n,\sigma^2)$ is the same as the Gaussian Symplectic
Ensemble from [Me].  

Whenever $\{X_j | j\in J\}$ is a family of random matrices in $\GSE(n,
\sigma^2)\cup \GSES(n, \sigma^2)$, we shall say that the $X_j$'s are
stochastically independent, if and only if the $V_j$'s, the $W_j$'s, the $Y_j$'s and the $Z_j$'s form a set of stochastically independent
random matrices.

\vspace{.2cm}

Finally, for the sake of completeness we define: \begin{itemize}
\item[(viii)] $\GRM^\H(n,\sigma^2)$ is the set of random matrices
  $Y\rightarrow M_n(\H)$ satisfying that
\[
Y = \unit\otimes V + j\otimes W + k\otimes X + l\otimes Z 
\]
for some stochastically independent random matrices $V, W, X$ and $Z$ from
$\GRMR(n,\frac{\sigma^2}{4})$. Equivalently,
\begin{equation}\label{GRMH}
Y = \frac{1}{\sqrt 2}(X_1 + \i X_2),
\end{equation}
for some stochastically independent random matrices
$X_1\in\GSE(n,\sigma^2)$ and $X_2\in\GSES(n,\sigma^2)$.
\end{itemize}

We apply Theorem~A to obtain:

{\bf Theorem~B.}
Let $r$, $s \in \N_0$ with $r+s\geq 1$ , and for each $n\in \N$, let $X_1^{(n)}, \ldots, X_{r+s}^{(n)}$ be stochastically independent random matrices defined on $(\Omega, \CF, P)$ such that $ X_1^{(n)}, \ldots,  X_r^{(n)}\in \GSE(n, \frac 1n)$ and $ X_{r+1}^{(n)},\ldots, X_{r+s}^{(n)} \in \GSES(n, \frac 1n)$. Furthermore, let $(\CA, \tau)$ be a $C\cc$-probability space with $\tau$ a faithful state on $\CA$, and let $\{x_1, \ldots, x_{r+s}\}$ be a semicircular system in $(\CA, \tau)$. Then there is a $P$-null set $N\subseteq\Omega$ such that for any $\omega\in \Omega\setminus N$ and for any polynomial $p$ in $r+s$ non-commuting variables:
\[
\lim_{n\rightarrow \infty}\|p(X_1^{(n)}(\omega), \ldots, X_{r+s}^{(n)}(\omega))\| = \|p(x_1, \ldots, x_{r+s})\|.
\]

\vspace{.2cm}

Haagerup and Thorbj\o rnsen apply [La, Proposition~4.1] and their main Theorem to prove that for $p\in \N$ and $Y_n\in \GRM(n, \frac 1n)$,
\[
\lim_{n\rightarrow\infty}\|Y_n^p\| =
\Bigg(\frac{(p+1)^{(p+1)}}{p^p}\Bigg)^{\frac 12}
\]
almost surely (cf. [HT2, Corollary~9.7]). Similarly, with the aid of [La, Proposition~4.1] and the identities (\ref{GRMR}) and (\ref{GRMH}), the following corollary follows from Theorem~A and Theorem~B:

{\bf Corollary.}
Let $p\in\N$, and for each $n\in\N$, let  $Y_n\in \GRMR(n, \frac 1n)$
($Y_n\in \GRM^\H (n, \frac 1n)$, respectively). Then, in both cases, 
  \[
  \lim_{n\rightarrow\infty}\|Y_n^p\| =
  \Bigg(\frac{(p+1)^{(p+1)}}{p^p}\Bigg)^{\frac 12}
  \]
  almost surely.

\section{Master equation and master inequality in the real case.}

Throughout this section let $r,\; m \in \mathbb{N}$ and $a_0, a_1, \ldots , a_r \in M_m(\mathbb{C})$ with $a_0^* = a_0$, and for each positive integer, $n$, let $Y_1^{(n)}, \ldots , Y_r^{(n)}$ be independent random matrices from $\textrm{GRM}^\mathbb{R}(n, \frac{1}{n})$. Then define the self-adjoint random matrix $S_n$ taking values in $M_m(\C)\otimes M_n(\C)$ by
\begin{equation}
  S_n = a_0\otimes \mathbf{1}_n + \sum_{j=1}^r (a_j\otimes Y_j^{(n)}+ a_j^*\otimes {Y_j^{(n)}}^*).
\end{equation}
Also, for $\lambda \in M_m(\C)$ set 
\[
\im \lambda = \frac{1}{2\i}(\lambda - \lambda\cc),
\]
and for $\im\lambda > 0$ (i.e. $\im \lambda$ is positive definite) define
\begin{equation}
  H_n(\lambda)  =  (\id_m\otimes \tr_n)[(\lambda\otimes\unit_n - S_n)^{-1}].
\end{equation}
We denote by $e_{kl}^{(m)}$ the matrix unit in $M_m(\C)$ with 1 at entry $(k,l)$ and 0 at all other entries. 

For an element $u \in M_m(\C)\otimes M_n(\C)$ we let $u^t$ denote the transpose of $u$ w.r.t. the natural identification of $M_m(\C)\otimes M_n(\C)$ with $M_{mn}(\C)$. Note that for $a\in M_m(\C)$ and $b\in M_n(\C)$, $(a\otimes b)^t = a^t\otimes b^t$. 

Finally, for any invertible square matrix $a$ we denote by $a^{-t}$ the transpose of $a^{-1}$. 

Applying the methods of Haagerup and Thorbj\o rnsen from \cite{HT2} we prove

\begin{mastereq}\label{MEQ}
  For every positive integer $n$ and every $\lambda \in M_m(\C)$ with $\im\lambda > 0$
  \begin{equation}
    \E\{(a_0-\lambda)H_n(\lambda) + \sum_{j=1}^r[a_j H_n(\lambda)a_j^* H_n(\lambda) + a_j^*H_n(\lambda)a_j H_n(\lambda)] + \mathbf{1}_m\} = -\frac{1}{n}R_n(\lambda),
  \end{equation}
  where
  \begin{equation}\label{MEQa}
    \begin{split}
    R_n(\lambda) =& \sum_{j=1}^r\sum_{k,l=1}^m a_j e_{kl}^{(m)} \E\{(\id_m\otimes \tr_n)[(\lambda\otimes\unit_n - S_n)^{-t}(e_{kl}^{(m)} a_j\otimes \unit_n)(\lambda\otimes\unit_n - S_n)^{-1}]\}\\
    &+ \sum_{j=1}^r\sum_{k,l=1}^m a_j^* e_{kl}^{(m)} \E\{(\id_m\otimes \tr_n)[(\lambda\otimes\unit_n - S_n)^{-t}(e_{kl}^{(m)} a_j^*\otimes \unit_n)(\lambda\otimes\unit_n - S_n)^{-1}]\}.
    \end{split}
  \end{equation}
\end{mastereq}
\proof
Consider a fixed $n \in \N$. For $1\leq j \leq r$ let $\{(Y_j^{(n)})_{kl}\}_{1\leq k,l\leq n}$ denote the entries of $ Y_j^{(n)}$. By [HT2, Lemma 3.3] and [HT2, Lemma 3.4],
\begin{eqnarray*}
  \E\{(Y_j^{(n)})_{kl}(\lambda\otimes\unit_n - S_n)^{-1}\} & = & \frac{1}{n}\E\Big\{\diff(\lambda\otimes\unit_n - S_n-t(a_j\otimes e_{kl}^{(n)}+a_j^*\otimes e_{lk}^{(n)}))^{-1}\Big\}\\[.2cm]
  & = & \frac{1}{n}\E\{(\lambda\otimes\unit_n - S_n)^{-1}(a_j\otimes e_{kl}^{(n)}+a_j^*\otimes e_{lk}^{(n)})(\lambda\otimes\unit_n - S_n)^{-1}\}.
\end{eqnarray*}
Since 
\[
  Y_j^{(n)} = \sum_{k,l=1}^n (Y_j^{(n)})_{kl} e_{kl}^{(n)},
\]
it follows that
\begin{equation*}
\begin{split}
  \E\{(a_j\otimes Y_j^{(n)})&(\lambda\otimes\unit_n - S_n)^{-1}\}=\\[.2cm]
  &\frac{1}{n} \sum_{k,l=1}^n (a_j\otimes e_{kl}^{(n)})\E\{(\lambda\otimes\unit_n - S_n)^{-1}(a_j\otimes e_{kl}^{(n)}+a_j^*\otimes e_{lk}^{(n)})(\lambda\otimes\unit_n - S_n)^{-1}\},
\end{split}
\end{equation*}
and
\begin{equation*}
\begin{split}
  \E\{(a_j^*\otimes {Y_j^{(n)}}^*)&(\lambda\otimes\unit_n - S_n)^{-1}\} =\\[.2cm]
  &\;\;\; \frac{1}{n} \sum_{k,l=1}^n (a_j^*\otimes e_{lk}^{(n)})\E\{(\lambda\otimes\unit_n - S_n)^{-1}(a_j\otimes e_{kl}^{(n)}+a_j^*\otimes e_{lk}^{(n)})(\lambda\otimes\unit_n - S_n)^{-1}\}.
\end{split}
\end{equation*}
As in the proof of [HT2, Lemma 3.5], for every $m\times m$ matrix $b$,
\[
  \sum_{k,l=1}^n (\unit_m\otimes e_{kl}^{(n)})(\lambda\otimes\unit_n - S_n)^{-1}(b\otimes e_{lk}^{(n)}) = n\;H_n(\lambda)b\otimes \unit_n. 
\]
Also observe that for any elementary tensor $x\otimes y \in M_m(\C)\otimes M_n(\C)$,
\begin{eqnarray*}
\sum_{k,l=1}^n (\unit_m\otimes e_{kl}^{(n)})(x\otimes y)(\unit_m\otimes e_{kl}^{(n)}) & = & x\otimes y^t = (x^t\otimes y)^t\\
& = &  \Big(\sum_{k,l=1}^m(e_{kl}^{(m)}\otimes \unit_n)(x\otimes y)(e_{kl}^{(m)}\otimes \unit_n)\Big)^t\\
& = & \sum_{k,l=1}^m(e_{lk}^{(m)}\otimes \unit_n)(x\otimes y)^t(e_{lk}^{(m)}\otimes \unit_n).
\end{eqnarray*}
Hence,
\[
  \sum_{k,l=1}^n (\unit_m\otimes e_{kl}^{(n)})(\lambda\otimes\unit_n - S_n)^{-1}(\unit_m\otimes e_{kl}^{(n)}) 
= \sum_{k,l=1}^m(e_{kl}^{(m)}\otimes \unit_n)(\lambda\otimes\unit_n - S_n)^{-t}(e_{kl}^{(m)}\otimes \unit_n).
\]
Combining the above observations we find that
\begin{equation*}
\begin{split}
  \E\{(a_j\otimes Y_j^{(n)})&(\lambda\otimes\unit_n - S_n)^{-1}\} = 
  \E\{(a_jH_n(\lambda)a_j^* \otimes \unit_n)(\lambda\otimes\unit_n - S_n)^{-1}\} + \\[.2cm]
  &\frac{1}{n}\sum_{k,l=1}^m(a_j e_{kl}^{(m)}\otimes \unit_n)\E\{(\lambda\otimes\unit_n - S_n)^{-t}(e_{kl}^{(m)}a_j\otimes \unit_n)(\lambda\otimes\mathbf{1}_n - S_n)^{-1}\},
\end{split}
\end{equation*}
and
\begin{equation*}
\begin{split}
  \E\{(a_j^*\otimes {Y_j^{(n)}}^*)&(\lambda\otimes\unit_n - S_n)^{-1}\} = 
 \E\{(a_j^*H_n(\lambda)a_j\otimes \unit_n)(\lambda\otimes\unit_n - S_n)^{-1}\} +\\[.2cm]
 &\frac{1}{n}\sum_{k,l=1}^m(a_j^* e_{kl}^{(m)}\otimes \unit_n)\E\{(\lambda\otimes\unit_n - S_n)^{-t}(e_{kl}^{(m)}a_j^*\otimes \unit_n)(\lambda\otimes\unit_n - S_n)^{-1}\}.
\end{split}
\end{equation*}

Applying now $(\id_m\otimes \tr_n)$ to both sides of the equations above we get
\begin{equation*}
\begin{split}
  \E\{(\id_m\otimes \tr_n)&[(a_j\otimes Y_j^{(n)})(\lambda\otimes\unit_n - S_n)^{-1}]\} = \E\{(a_jH_n(\lambda)a_j^*H_n(\lambda)\} +\\[.2cm]
  &\frac{1}{n}\sum_{k,l=1}^ma_je_{kl}^{(m)}\E\{(\id_m\otimes \tr_n)[(\lambda\otimes\unit_n - S_n)^{-t}(e_{kl}^{(m)}a_j\otimes \unit_n)(\lambda\otimes\mathbf{1}_n - S_n)^{-1}]\},
\end{split}
\end{equation*}
\begin{equation*}
\begin{split}
  \mathbb{E}\{(\textrm{id}_m\otimes \textrm{tr}_n)&[(a_j^*\otimes {Y_j^{(n)}}^*)(\lambda\otimes\mathbf{1}_n - S_n)^{-1}]\} = \mathbb{E}\{(a_j^*H_n(\lambda)a_jH_n(\lambda)\} +\\[.2cm] 
  &\frac{1}{n}\sum_{k,l=1}^m(a_j^*e_{kl}^{(m)}\mathbb{E}\{(\textrm{id}_m\otimes \textrm{tr}_n)[(\lambda\otimes\mathbf{1}_n - S_n)^{-t}(e_{kl}^{(m)}a_j^*\otimes \mathbf{1}_n)(\lambda\otimes\mathbf{1}_n - S_n)^{-1}]\}.
\end{split}
\end{equation*}
Since
\begin{equation*}
  \mathbb{E}\{(a_0-\lambda)H_n(\lambda)\} =  \mathbb{E}\{(\textrm{id}_m\otimes \textrm{tr}_n)((a_0-\lambda)\otimes \mathbf{1}_n)(\lambda\otimes\mathbf{1}_n - S_n)^{-1}\},
\end{equation*}
we conclude that
\begin{equation*}
\begin{split}
  \mathbb{E}\{(a_0&-\lambda)H_n(\lambda)+\sum_{j=1}^ra_jH_n(\lambda)a_j^*H_n(\lambda)+a_j^*H_n(\lambda)a_jH_n(\lambda)\} =\\
& \mathbb{E}\{(\textrm{id}_m\otimes \textrm{tr}_n)[((a_0-\lambda)\otimes \mathbf{1}_n + \sum_{j=1}^ra_j\otimes Y_j^{(n)} +a_j^*\otimes {Y_j^{(n)}}^*) (\lambda\otimes\mathbf{1}_n - S_n)^{-1}]\} - \frac{1}{n}R_n(\lambda).
\end{split}
\end{equation*}
Hence,
\begin{equation*}
  \mathbb{E}\{(a_0-\lambda)H_n(\lambda) + \sum_{j=1}^r[a_j H_n(\lambda)a_j^* H_n(\lambda) + a_j^*H_n(\lambda)a_j H_n(\lambda)] +  \mathbf{1}_m\}  = - \frac{1}{n}R_n(\lambda),
\end{equation*}
and the proof is complete.$\endproof$

\vspace{.2cm}

Now, let $(\CA, \tau)$ be a $C^*$-probability space, and let $y_1, \ldots, y_r$ be a circular system in $(\CA, \tau)$, i.e. with
\begin{eqnarray*}
  x_j = &\frac{1}{\sqrt{ 2}}(y_j + y_j^*), &\;\;\;(1\leq j \leq r)\\
  x_{j +r} =   &\frac{1}{\i\sqrt{ 2}}(y_j - y_j^*), &\;\;\;(1\leq j \leq r)
\end{eqnarray*}
$x_1, \ldots, x_{2r}$ form a semicircular system. Define
\begin{equation}
   s = a_0\otimes \unit_\CA + \sum_{j=1}^r (a_j\otimes y_j + a_j^*\otimes y_j^*). 
\end{equation}

Note that
\[ 
   s = a_0\otimes \unit_\CA + \sum_{j=1}^r\Big( \frac{1}{\sqrt{ 2}}(a_j+a_j^*)\otimes x_j + \frac{1}{\i\sqrt{ 2}}(a_j-a_j^*)\otimes x_{j+r}\Big).
\]
For $\lambda \in M_m(\C)$ with $\im\lambda > 0$ define
\begin{equation}
  G(\lambda) = (\id_m\otimes \tau)[(\lambda\otimes\unit_\CA - s)^{-1}].
\end{equation}
According to [HT2, Lemma 5.4] we have the following identity
\begin{equation}
  \lambda = a_0 + \sum_{j=1}^r\Big[\frac{1}{\sqrt{ 2}}(a_j+a_j^*)G(\lambda)\frac{1}{\sqrt{ 2}}(a_j+a_j^*) + \frac{1}{\i\sqrt{ 2}}(a_j-a_j^*)G(\lambda)\frac{\textrm{i}}{\sqrt{ 2}}(a_j-a_j^*)\Big] + G(\lambda)^{-1},
\end{equation}
or equivalently,
\begin{equation}\label{GEQ}
  \lambda = a_0 + \sum_{j=1}^r[a_jG(\lambda)a_j^* + a_j^*G(\lambda)a_j] + G(\lambda)^{-1}.
\end{equation}

\vspace{.2cm}

\begin{remark}\label{psi}
Let $r$ and $n$ be positive integers, and define an isomorphism $\psi_0 : M_n(\C)\rightarrow \C^{n^2}$ by
\[
\psi_0((a_{kl})_{1\leq k,l\leq n}) = ((a_{kl})_{1\leq k,l\leq n}).
\]
for $(a_{kl})_{1\leq k,l\leq n}\in M_n(\C)$. $\psi_0$ has a natural extension to a linear isomorphism between $M_n(\C)^r$ and $\C^{rn^2}$, which we denote by $\psi$:
\[
\psi(A_1, \ldots, A_r) = (\psi_0(A_1), \ldots, \psi_0(A_r)), \;\;\;\;\;\;(A_1, \ldots , A_r \in M_n(\C)).
\]

Define a norm $\|\cdot\|_e$ on $M_n(\C)^r$ by
\[
\|(A_1, \ldots, A_r)\|_e^2 = \Tr_n\Big(\sum_{i=1}^r A_i\cc A_i\Big),  \;\;\;\;\;\;(A_1, \ldots , A_r \in M_n(\C)),
\]
and note that $\psi$ is an isometry with respect to this norm and the Euclidian norm on  $\C^{rn^2}$. 

If $Y_1^{(n)}, \ldots, Y_r^{(n)}$ are independent random matrices from
$\GRMR(n,\frac 1n)$, then $\Y = \psi(Y_1^{(n)}, \ldots, Y_r^{(n)})$ is a
random variable taking values in $\R^{rn^2}$, and the distribution of $\Y$
on $\R^{rn^2}$ is $\nu\otimes \cdots \otimes \nu$ ($rn^2$ terms), where
$\nu$ denotes the Gaussian distribution $N(0, \frac 1n)$ on $\R$. By the
Gaussian Poincar\'e inequality (cf. [HT2, Remark~4.3]), if $\tilde{f}: \R^{rn^2} \rightarrow \C$ is a $C^1$-function such that $\tilde{f}$ and $\grad\tilde{f}$ are polynomially bounded, then, with $f = \tilde{f}(\psi): M_n(\C)^r \rightarrow \C$:
\[
\V\{f(Y_1^{(n)}, \ldots, Y_r^{(n)})\} \leq \frac 1n \E\{\|\grad f (Y_1^{(n)}, \ldots, Y_r^{(n)})\|_e^2\},
\] 
where $\V\{g\} = \E\{|g - \E\{g\}|^2\}$ for any function $g\in L^2(\Omega, P)$.
\end{remark}

\vspace{.2cm}

\begin{lemma}\label{sumnorm}
Let $r$, $m$ and $n$ be positive integers, let $a_1, \ldots, a_r \in M_m(\C)$, and let $w_1, \ldots, w_r \in M_n(\C)$. Then, with $w = (w_1, \ldots, w_r)$,
\[
 \Big\|\sum_{i=1}^r a_i\otimes w_i+ a_i^*\otimes w_i^*\Big\|_{2, \Tr_m\otimes\Tr_n} \leq 2 m^{1/2}\Big(\sum_{i=1}^r\|a_i\|^2\Big)^{1/2} \|w\|_e,
\]
where $\|\cdot \|_{2,\Tr_m\otimes \Tr_n}$ denotes the Hilbert-Schmidt norm
on $M_m(\C)\otimes M_n(\C)$.
\end{lemma}

\proof
This is a simple application of the Cauchy-Schwartz inequality for the standard inner product on $\C^r$:
\begin{eqnarray*}
  \Big\|\sum_{i=1}^r a_i\otimes w_i+ a_i^*\otimes w_i^*\Big\|_{2, \Tr_m\otimes\Tr_n} & \leq & 2 \sum_{i=1}^r\|a_i\|_{2, \Tr_m}\|w_i\|_{2, \Tr_n}\\
& \leq & 2 \Big(\sum_{i=1}^r\|a_i\|_{2, \textrm{Tr}_m}^2\Big)^{1/2} \Big(\sum_{i=1}^r\|w_i\|_{2, \textrm{Tr}_n}^2\Big)^{1/2}\\
& = & 2 \Big(\sum_{i=1}^r \textrm{Tr}_m(a_i^*a_i)\Big)^{1/2}\|w\|_e\\
& \leq &  2 m^{1/2}\Big\|\sum_{i=1}^ra_i\cc a_i\Big\|^{1/2}  \|w\|_e\\
& \leq &  2 m^{1/2}\Big(\sum_{i=1}^r\|a_i\|^2\Big)^{1/2} \|w\|_e.\;\;\;\;\;\;\endproof
\end{eqnarray*}

\vspace{.2cm}

\begin{masterineq}\label{MIEQ}
  There is a constant $C_1 <\infty$ such that for every positive integer $n$ and for every $\lambda \in M_m(\mathbb{C})$ with $\textrm{Im}\lambda > 0$,
\begin{equation}
\begin{split}
  \Big\|a_0 + \sum_{j=1}^r[a_jG_n(\lambda)a_j^* + a_j^*G_n(\lambda)a_j] + G_n(\lambda)^{-1} - \lambda &+ \frac{1}{n}R_n(\lambda) G_n(\lambda)^{-1}\Big\|\\[.2cm]
  &\leq \frac{C_1}{n^2}(\|\lambda\|+K)^2\|(\textrm{Im}\lambda)^{-1}\|^5,
\end{split}
\end{equation}
where $K = \|a_0\|+ 16\sum_{j=1}^r\|a_j\|$.
\end{masterineq}

\proof
Let $n \in \N$, and let $\lambda \in M_m(\C)$ with $\im\lambda > 0$. Define
\begin{equation}
  K_n(\lambda) = H_n(\lambda) - G_n(\lambda).
\end{equation}
According to  Theorem~\ref{MEQ}
\begin{equation*}
\begin{split}
  \E\Big\{\sum_{j=1}^r&a_jK_n(\lambda)a_j^*K_n(\lambda) + a_j^*K_n(\lambda)a_jK_n(\lambda)\Big\} = \\
&- \Big(\sum_{j=1}^r[a_jG_n(\lambda)a_j^*G_n(\lambda) + a_j^*G_n(\lambda)a_jG_n(\lambda)] + (a_0 -\lambda)G_n(\lambda) + \mathbf{1}_m\Big) - \frac{1}{n}R_n(\lambda),
\end{split}
\end{equation*}
i.e.
\begin{equation*}
\begin{split}
 -\E\Big\{\sum_{j=1}^r&a_jK_n(\lambda)a_j^*K_n(\lambda) + a_j^*K_n(\lambda)a_jK_n(\lambda)\Big\}G_n(\lambda)^{-1} = \\
&a_0 + \sum_{j=1}^r[a_jG_n(\lambda)a_j^* + a_j^*G_n(\lambda)a_j] + G_n(\lambda)^{-1} -\lambda + \frac{1}{n}R_n(\lambda)G_n(\lambda)^{-1},
\end{split}
\end{equation*}
which implies that
\begin{equation*}
\begin{split}
  \Big\|a_0 + \sum_{j=1}^r[a_jG_n(\lambda)a_j^* +& a_j^*G_n(\lambda)a_j] + G_n(\lambda)^{-1} - \lambda + \frac{1}{n}R_n(\lambda) G_n(\lambda)^{-1}\Big\|\\
  &\leq \E\Big\{\Big\|\sum_{j=1}^ra_jK_n(\lambda)a_j^* + a_j^*K_n(\lambda)a_j\Big\|\|K_n(\lambda)\|\Big\}\|G_n(\lambda)^{-1}\|.
\end{split}
\end{equation*}
The mappings $v\mapsto a_jva_j^*$, $v\mapsto a_j^*va_j$, $1\leq j\leq r$,
are all completely positive. Hence, the sum of these mappings is also
completely positive, so it attains it norm at the unit of $M_m(\C)$, and we have
\begin{equation*}
\begin{split}
  \Big\|a_0 + \sum_{j=1}^r[a_jG_n(\lambda)a_j^* + a_j^*G_n(\lambda)a_j] +& G_n(\lambda)^{-1} - \lambda + \frac{1}{n}R_n(\lambda) G_n(\lambda)^{-1}\Big\|\\
  &\leq\Big\|\sum_{j=1}^ra_ja_j^* + a_j^*a_j\Big\|\;\E\{\|K_n(\lambda)\|^2\}\;\|G_n(\lambda)^{-1}\| \\
  &\leq 2 \sum_{j=1}^r \|a_j\|^2 \;\E\{\|K_n(\lambda)\|_{2, \Tr_m}^2\}\;\|G_n(\lambda)^{-1}\|,
\end{split}
\end{equation*}
where $\|K_n(\lambda)\|_{2, \Tr_m}$ denotes the Hilbert-Schmidt norm of $K_n(\lambda)$.

Now, for $Y\in \GRMR(n,\frac{1}{n})$ we can choose $Z \in \GRM(n,\frac{1}{n})$ such that 
\[
  Y=\frac{1}{\sqrt 2}(Z+\overline{Z}),
\]
where $\overline{Z} = (Z\cc)^t$.
The random matrix $Z$ may be expressed in terms of two indenpendent random matrices $X_1, X_2 \in \SGRM(n,\frac{1}{n})$:
\[
  Z = \frac{1}{\sqrt{2}}(X_1+\textrm{i}X_2).
\]
It follows from [HT2, Lemma 5.1] that
\[
  \E\{\|Y\|\}\leq 2\E\{\|X_1\|\} \leq 8.
\]

By arguments similar to those presented in [HT2, Proof of Proposition 5.2],
\begin{equation}\label{Gn-1}
  \|G_n(\lambda)^{-1}\|\leq (\|\lambda\| + K)^2\|(\im\lambda)^{-1}\|,
\end{equation}
where $K = \|a_0\|+ 16\sum_{j=1}^r\|a_j\|$. 

To get an estimate of $\E\{\|K_n(\lambda)\|_{2, \Tr_m}^2\}$ we follow the lines of [HT2, Proof of Theorem 4.5] and obtain:  
\begin{equation}\label{K}
  \E\{\|K_n(\lambda)\|_{2, \Tr_m}^2\} \leq \frac{1}{n}\sum_{j,k =1}^m \E\{\|(\grad f_{n,j,k})(Y_1^{(n)}, \ldots, Y_r^{(n)})\|_e^2\},
\end{equation}
where each $f_{n,j,k}: M_n(\C)^r\rightarrow \C$ is defined by
\[
  f_{n,j,k}(v_1, \ldots, v_r) = m(\tr_m\otimes\tr_n)[(e_{jk}^{(m)}\otimes \unit_n)(\lambda \otimes \unit_n - a_0\otimes \mathbf{1}_n - \sum_{i=1}^r a_i\otimes v_i+ a_i^*\otimes v_i^*)^{-1}].
\]
Let $1\leq j, k\leq m$, let $v= (v_1, \ldots , v_r) \in  M_n(\mathbb{C})^r$, and let $w= (w_1, \ldots , w_r) \in  M_n(\mathbb{C})^r$ with $\|w\|_e = 1$. Proceeding as in [HT2, Proof of Theorem 4.5] we find that
\[
  \Big|\diff f_{n,j,k}(v+tw)\Big|^2 \leq  \frac{1}{n}\Big\|\sum_{i=1}^r a_i\otimes w_i+ a_i^*\otimes w_i^*\Big\|_{2, \textrm{Tr}_m\otimes\textrm{Tr}_n}^2\|(\textrm{Im}\lambda)^{-1}\|^4.
\]
Hence, by  Lemma~\ref{sumnorm}
\begin{equation}
  \Big|\diff f_{n,j,k}(v+tw)\Big |^2 \leq  \frac{4m}{n}   \sum_{i=1}^r \|a_i\|^2  \|(\textrm{Im}\lambda)^{-1}\|^4,
\end{equation}
and since $w$ was arbitrary,
\begin{equation}
  \|\grad f_{n,j,k}(v)\|_e^2 \leq \frac{4m}{n} \sum_{i=1}^r \|a_i\|^2   \|(\im\lambda)^{-1}\|^4.
\end{equation}
Inserting this into (\ref{K}) we obtain
\begin{equation}
  \mathbb{E}\{\|K_n(\lambda)\|_{2, \textrm{Tr}_m}^2\} \leq \frac{4m^3}{n^2}  \sum_{i=1}^r \|a_i\|^2
  \|(\textrm{Im}\lambda)^{-1}\|^4.
\end{equation}
Finally,
\begin{equation*}
\begin{split}
  \Big\|a_0 + \sum_{j=1}^r&[a_jG_n(\lambda)a_j^* + a_j^*G_n(\lambda)a_j] + G_n(\lambda)^{-1} - \lambda + \frac{1}{n}R_n(\lambda) G_n(\lambda)^{-1}\Big\|\\
  &\leq 2 \sum_{j=1}^r\|a_j\|^2 \;  \frac{4m^3}{n^2} \; \sum_{i=1}^r
  \|a_i\|^2  \;  \|(\textrm{Im}\lambda)^{-1}\|^4  (\|\lambda\| + K)^2\|(\textrm{Im}\lambda)^{-1}\| \\
 &=\frac{8m^3}{n^2} \Big(\sum_{j=1}^r \|a_j\|^2\Big)^2  (\|\lambda\| + K)^2 \|(\textrm{Im}\lambda)^{-1}\|^5,
\end{split}
\end{equation*}
from which the theorem follows.$\endproof$

\begin{cor}\label{CMIEQ}
   There is a constant $C_1' <\infty$ such that for every positive integer $n$ and for every $\lambda \in M_m(\mathbb{C})$ with $\textrm{Im}\lambda > 0$,
\begin{equation}
\begin{split}
  \Big\|a_0 + \sum_{j=1}^r[a_jG_n(\lambda)a_j^* + a_j^*G_n(\lambda)a_j] &+ G_n(\lambda)^{-1} - \lambda\Big\|\\
  &\leq \frac{C_1'}{n}(\|\lambda\|+K)^2(\|(\textrm{Im}\lambda)^{-1}\|^5+\|(\textrm{Im}\lambda)^{-1}\|^3),
\end{split}
\end{equation}
where $K = \|a_0\|+ 16\sum_{j=1}^r\|a_j\|$.
\end{cor}

\proof
By application of [HT2, Lemma 3.1] we find that
\begin{eqnarray*}
  \|R_n(\lambda)\|&  \leq & 2\sum_{k,l =1}^m\sum_{j=1}^r\|a_j\|\mathbb{E}\{\|(\lambda\otimes\mathbf{1}_n-S_n)^{-t}\|\|a_j\|\|(\lambda\otimes\mathbf{1}_n-S_n)^{-1}\|\}\\
& \leq & 2m^2\Big(\sum_{j=1}^r \|a_j\|^2\Big)\|(\textrm{Im}\lambda)^{-1}\|^2.
\end{eqnarray*}

Then, by Theorem~\ref{MIEQ} and by (\ref{Gn-1})
\begin{equation*}
\begin{split}
\Big\|a_0 + \sum_{j=1}^r&[a_jG_n(\lambda)a_j^* + a_j^*G_n(\lambda)a_j] + G_n(\lambda)^{-1} - \lambda\Big\|\\
&\leq\frac{C_1}{n^2}(\|\lambda\|+K)^2\|(\textrm{Im}\lambda)^{-1}\|^5 + \frac{1}{n}\|R_n(\lambda)\|\|G_n(\lambda)^{-1}\| \\
&\leq\frac{C_1}{n}(\|\lambda\|+K)^2\|(\textrm{Im}\lambda)^{-1}\|^5 + \frac{2m^2}{n}\Big(\sum_{j=1}^r \|a_j\|^2\Big)\|(\textrm{Im}\lambda)^{-1}\|^2(\|\lambda\|+K)^2\|(\textrm{Im}\lambda)^{-1}\|,
\end{split}
\end{equation*}
and the corollary follows.$\endproof$

\section{Estimation of $\|G_n(\lambda)-G(\lambda)\|$.}

We stick to the notation introduced in the previous section and define the subset $\CO$ of $M_m(\C)$ by
\begin{equation}
  \CO = \{\lambda \in M_m(\C)| \im\lambda > 0\},
\end{equation}
and for $\lambda \in \CO$ put
\begin{eqnarray}
  \Lambda_n(\lambda) & = & a_0 + \sum_{j=1}^r[a_jG_n(\lambda)a_j^* + a_j^*G_n(\lambda)a_j] + G_n(\lambda)^{-1},\\
  \eps(\lambda) & = & \|(\im\lambda)^{-1}\|^{-1}.
\end{eqnarray}
Finally let
\begin{eqnarray}
  \CO_n'& =& \Big\{\lambda \in \CO\Big| \frac{C_1'}{n}(K+\|\lambda\|)^2(\eps(\lambda)^{-6}+\eps(\lambda)^{-4})<\frac{1}{2}\Big\}.
\end{eqnarray}

By application of Corollary~\ref{CMIEQ} and the methods of [HT2, Proof of Lemma~5.5] one finds that for any $\lambda \in \CO_n'$,
\[
  \im\Lambda_n(\lambda) \geq \frac{\eps(\lambda)}{2}\unit_m.
\]

Thus, $\Lambda_n(\lambda) \in \CO$ with
\begin{equation}\label{Lambda_n}
\|(\im\Lambda_n(\lambda))^{-1}\| \leq 2\|(\im\lambda)^{-1}\|, 
\end{equation}
and by (\ref{GEQ}), 
\begin{equation*}
\begin{split}
  a_0 + \sum_{j=1}^r[a_jG_n(\lambda)a_j^* + a_j^*G_n(\lambda)a_j] &+ G_n(\lambda)^{-1} = \\
  & a_0 + \sum_{j=1}^r[a_jG(\Lambda_n(\lambda))a_j^* + a_j^*G(\Lambda_n(\lambda))a_j] + G(\Lambda_n(\lambda))^{-1}.
\end{split}
\end{equation*}

\vspace{.2cm}

\begin{prop}\label{biglambda}
Let $n\in \N$. Then for all $\lambda \in \CO_n'$,
\begin{equation*}
  G(\Lambda_n(\lambda))=G_n(\lambda).
\end{equation*}
\end{prop}

\proof The proof is almost identical to [HT2, Proof of Proposition~5.5]. Only a few modifications are necessary to make it work in this case too. 

Making use of the fact that the function $t \mapsto (K+t)^2(t^{-6}+t^{-4})$, $t>0$, is continuous and strictly decreasing, it follows as in \cite{HT2} that 
\begin{itemize}
  \item[(a)] $\CO_n'$ is an open connected subset of $M_m(\C)$.
\end{itemize}
With 
\[
  \CO_n'' = \{\lambda \in \CO_n'| \eps(\lambda)>2(\sum_{j=1}^r\|a_j\|^2)^{\frac 12}\},
\]
which is an open, nonempty subset of $ \CO_n'$, one gets, as in \cite{HT2}, that
\begin{itemize}
  \item[(b)] $G(\Lambda_n(\lambda))=G_n(\lambda)$ for all $\lambda \in  \CO_n''$.
\end{itemize}
Finally, apply the principle of uniqueness of analytic continuation.$\endproof$

\vspace{.2cm}

Taking Corollary~\ref{CMIEQ} into account and proceeding as in [HT2, Proof of Theorem~5.7] one gets:

\begin{thm}\label{Gn-G}
There is a constant $C_2 <\infty$ such that for any $\lambda \in \CO$ and for any positive integer $n$
  \begin{equation}
    \|G_n(\lambda)-G(\lambda)\|\leq \frac{C_2}{n}(\|\lambda\|+K)^2(\|(\im\lambda)^{-1}\|^7 +\|(\im\lambda)^{-1}\|^5),
  \end{equation}
where $K = \|a_0\|+ 16\sum_{j=1}^r\|a_j\|$.
\end{thm}

\vspace{.2cm}

Before stating the next corollary we introduce some notation. Let $\CU$ be
an open subset of $M_m(\C)$, and let $\phi: \CU \rightarrow M_m(\C)$ be a
complex differentiable map. For $v \in \CU$ we denote by $\phi'(v)$ the
differential of $\phi$ at $v$, i.e. $\phi'(v) : M_m(\C)
\rightarrow M_m(\C)$ is the (unique) linear map satisfying that for every
differentiable curve $\alpha$ defined in a neighbourhood of zero, $\alpha
: (-\eps, \eps) \rightarrow \CU$, with
$\alpha(0) = v$, the tangent vector $\diff \phi(\alpha(t))$ is given by
\[
\diff \phi(\alpha(t)) = \phi'(v)[\alpha'(0)].
\]

\vspace{.2cm}

\begin{cor}\label{diff}
Let $\lambda \in \CO$, and let $C_2$ be as in Theorem~\ref{Gn-G}. Then for every positive integer $n$ we have
\begin{equation}
\|G_n'(\lambda)-G'(\lambda)\| \leq \frac{576 C_2}{n}(\|\lambda\|+K)^2
(\|(\im\lambda)^{-1}\|^8 +\|(\im\lambda)^{-1}\|^6).
\end{equation}
\end{cor}

\proof
Let $n\in\N$ and let $x\in M_m(\mathbb{C})$ with $\|x\|=1$. Note that for
every complex number $z$
\[
\im(\lambda+zx) = \im\lambda + \im(zx) \geq \im\lambda -|z|\unit_m.
\]
So $\lambda + zx \in \CO$ if $|z| < \eps(\lambda) =
\|(\im\lambda)^{-1}\|^{-1}$. Hence, the map
\[
z \mapsto G_n(\lambda +zx)-G(\lambda + zx)
\]
is well-defined (and analytic) in the open disc in $\C$ of radius
$\|(\im\lambda)^{-1}\|^{-1}$ centered at zero. Put $r = \frac 12
\|(\im\lambda)^{-1}\|^{-1}$, and define a path $\gamma$ by
$\gamma(t) = r\e^{\i t}$, $t \in [0, 2\pi]$. Then, according to the Cauchy
formulas for vector valued analytic functions,
\begin{equation*}
\begin{split}
\|(G_n'(\lambda)-&G'(\lambda))[x]\|\\
 & =  \Big\|\diffz(G_n(\lambda
+zx)-G(\lambda +zx))\Big\|\\
&= \Big\|\frac{1}{2\pi \i}\int_\gamma \frac{G_n(\lambda + \zeta x)-G(\lambda +
  \zeta x)}{\zeta ^2}\d\zeta \Big\|\\
& \leq \frac{1}{2\pi}\int_0^{2\pi}\frac{\|G_n(\lambda +
  \gamma(t)x-G(\lambda + \gamma(t)x)\|}{r}\d t\\
& \leq  \frac 1r \max_{0\leq t \leq 2\pi}\{\|G_n(\lambda + \gamma(t)x)
-G(\lambda + \gamma(t)x)\|\}\\
& \leq  \frac{C_2}{nr}\max_{0\leq t \leq 2\pi}\{(\|\lambda
+\gamma(t)x\|+K)^2\cdot(\|(\im(\lambda
+\gamma(t)x))^{-1}\|^7 +\|(\im(\lambda+\gamma(t)x))^{-1}\|^5)\}.
\end{split}
\end{equation*}

Now, for $t \in [0, 2\pi]$
\begin{equation*}
  \im(\lambda + \gamma(t)x)  =  \im\lambda + r\im(e^{it}x)\geq (\eps(\lambda)-r)\unit_m = \frac 12 \|(\im\lambda)^{-1}\|^{-1}\unit_m.
\end{equation*}
Hence,
\[  
  \|(\im(\lambda +\gamma(t)x))^{-1}\| \leq  2 \|(\im\lambda)^{-1}\|.
\]
Also, since 
\[
\|\lambda\|\|(\im \lambda)^{-1}\| \geq \|\im\lambda\|\|(\im \lambda)^{-1}\|\geq 1,
\]
we have
\[
  \|\lambda + \gamma(t)x\| \leq \|\lambda\| + r = \|\lambda\| + \frac{1}{2\|(\textrm{Im}\lambda)^{-1}\|}\leq \frac 32 \;\|\lambda\|.
\]
It follows that
\begin{eqnarray*}
  \|(G_n'(\lambda)-G'(\lambda))[x]\| & \leq & \frac{C_2}{n}\Big(\frac 32\|\lambda\|+K\Big)^2(2^8\|(\im\lambda)^{-1}\|^8 + 2^6\|(\im\lambda)^{-1}\|^6)\\
&\leq& \frac{256 C_2}{n}\frac 94 (\|\lambda\|+K)^2
(\|(\im\lambda)^{-1}\|^8 +\|(\im\lambda)^{-1}\|^6)\\
& = &\frac{576 C_2}{n}(\|\lambda\|+K)^2
(\|(\im\lambda)^{-1}\|^8 +\|(\im\lambda)^{-1}\|^6).\;\;\;\;\;\;\endproof
\end{eqnarray*}

\section{Estimation of $\|G(\lambda)-G_n(\lambda)+\frac 1n L(\lambda)\|$.}

For $a_0, a_1, \ldots, a_r$ in $M_m(\C)$ with $a_0 = a_0\cc$ (as in the previous sections) define
\[
\tilde{a_j} = \begin{pmatrix}
\overline{a_j} & 0\\
0 & a_j
\end{pmatrix} \in M_{2m}(\C),\;\;\;\;\;\; (0\leq j\leq r),
\]
where $\overline{a_j}$ is the matrix obtained by conjugation of the entries of $a_j$. Note that $\overline{a_j} = {(a_j\cc)}^t$. 

For $\lambda \in M_m(\mathbb{C})$ with $\textrm{Im}\lambda > 0$ set
\[
\tilde{\lambda} = \begin{pmatrix}
\lambda^t & 0\\
0 & \lambda
\end{pmatrix} \in M_{2m}(\C), 
\]
and for $Y_1^{(n)}, \ldots, Y_r^{(n)}$ stochastically independent random matrices from the class $\GRMR(n, \frac 1n)$ define
\begin{eqnarray}
\tilde{S}_n & = & \tilde{a_0}\otimes \mathbf{1}_n + \sum_{j=1}^r (\tilde{a_j}\otimes Y_j^{(n)}+ \tilde{a_j}^*\otimes {Y_j^{(n)}}^*),\\
\tilde{s} & = &  \tilde{a_0}\otimes \unit_\CA + \sum_{j=1}^r (\tilde{a_j}\otimes y_j + \tilde{a_j}^*\otimes y_j^*),\\
\tilde{G}_n(\tilde{\lambda}) & = & \E\{(\id_{2m}\otimes \tr_n)[(\tilde{\lambda}\otimes\unit_n - \tilde{S}_n)^{-1}]\},\\
\tilde{G}(\tilde{\lambda}) & = & (\id_{2m}\otimes \tau)[(\tilde{\lambda}\otimes\unit_\CA - \tilde{s})^{-1}].
\end{eqnarray}

Note that $\|\tilde{a_j}\| = \|a_j\|$, $0\leq j\leq r$,  $\|\tilde{\lambda}\| = \|\lambda\|$, and that $\im \tilde{\lambda} > 0$ with $\|(\im \tilde{\lambda})^{-1}\| = \|(\im \lambda)^{-1}\|$. Thus, it follows from the results obtained this far that there is a constant, $\tilde{C}_2< \infty$, such that for every positive integer $n$ and every $\lambda \in M_m(\mathbb{C})$ with $\textrm{Im}\lambda > 0$,
\begin{equation}
    \|\tilde{G}_n(\tilde{\lambda})-\tilde{G}(\tilde{\lambda})\|\leq \frac{\tilde{C}_2}{n}(\|\lambda\|+K)^2(\|(\im\lambda)^{-1}\|^7 +\|(\im\lambda)^{-1}\|^5),
  \end{equation}
where $K = \|a_0\|+ 16\sum_{j=1}^r\|a_j\|$. The proof of Corollary~\ref{diff} also carries over, so
\begin{equation}\label{diff1}
\|\tilde{G}_n'(\tilde{\lambda})-\tilde{G}'(\tilde{\lambda})\| \leq \frac{576 \tilde{C_2}}{n}(\|\lambda\|+K)^2
(\|(\im\lambda)^{-1}\|^8 +\|(\im\lambda)^{-1}\|^6).
\end{equation}

\vspace{.2cm}

\begin{lemma}\label{Rn}
For every positive integer $n$ and every $\lambda \in M_m(\C)$ with $\im \lambda > 0$,
\begin{eqnarray*}
R_n(\lambda) & = & - \sum_{j=1}^r\sum_{k,l=1}^ma_je_{kl}^{(m)}(\Tr_2\otimes\id_m)\Bigg\{\tilde{G}_n'(\tilde{\lambda})\Bigg[
\begin{pmatrix}
0 & e_{kl}^{(m)}a_j\\
0 & 0
\end{pmatrix}
\Bigg]
\begin{pmatrix}
0 & 0\\
\unit_m & 0
\end{pmatrix}
\Bigg\}\\
& &  - \sum_{j=1}^r\sum_{k,l=1}^ma_j\cc e_{kl}^{(m)}(\Tr_2\otimes\id_m)\Bigg\{\tilde{G}_n'(\tilde{\lambda})\Bigg[
\begin{pmatrix}
0 & e_{kl}^{(m)}a_j\cc\\
0 & 0
\end{pmatrix}
\Bigg]
\begin{pmatrix}
0 & 0\\
\unit_m & 0
\end{pmatrix}
\Bigg\}.
\end{eqnarray*}
\end{lemma}

\proof
By definition
\begin{equation*}
\begin{split}
    R_n(\lambda) = & \sum_{j=1}^r\sum_{k,l=1}^m a_j e_{kl}^{(m)} \mathbb{E}\{(\textrm{id}_m\otimes \textrm{tr}_n)[(\lambda\otimes\mathbf{1}_n - S_n)^{-t}(e_{kl}^{(m)} a_j\otimes \mathbf{1}_n)(\lambda\otimes\mathbf{1}_n - S_n)^{-1}]\}\\
    &+\sum_{j=1}^r\sum_{k,l=1}^m a_j^* e_{kl}^{(m)} \mathbb{E}\{(\textrm{id}_m\otimes \textrm{tr}_n)[(\lambda\otimes\mathbf{1}_n - S_n)^{-t}(e_{kl}^{(m)} a_j^*\otimes \mathbf{1}_n)(\lambda\otimes\mathbf{1}_n - S_n)^{-1}]\}.
\end{split}
\end{equation*}

Note that $(\lambda\otimes\unit_n-S_n)^{-t} = (\lambda^t\otimes\unit_n-S_n^t)^{-1}$. Moreover, since the random matrices $Y_j^{(n)}$ have real entries,
\[
S_n^t = a_0^t\otimes\unit_n + \sum_{j=1}^r\Big(a_j^t\otimes{Y_j^{(n)}}\cc + (a_j\cc)^t\otimes Y_j^{(n)}\Big).
\]
As $a_0 = a_0\cc$, it follows that
\[
S_n^t = \overline{a_0}\otimes\unit_n + \sum_{j=1}^r\Big(\overline{a_j}\otimes Y_j^{(n)} + \overline{a_j}^t\otimes {Y_j^{(n)}}\cc\Big).
\]
Hence,
\[
\tilde{S_n} = 
\begin{pmatrix}
S_n^t & 0\\
0 & S_n
\end{pmatrix}.
\]
Now, standard matrix manipulations reveal that for every $x \in M_m(\C)$,
\[
\begin{pmatrix}
0 & (\lambda\otimes\unit_n - S_n)^{-t}(x\otimes \unit_n)(\lambda\otimes\unit - S_n)^{-1}\\
0 & 0
\end{pmatrix} = 
(\tilde{\lambda}\otimes\unit_n - \tilde{S}_n)^{-1}(\tilde{x}\otimes\unit_n)(\tilde{\lambda}\otimes\unit_n - \tilde{S}_n)^{-1},
\]
where
\[
\tilde{x} =\begin{pmatrix}
0 & x\\
0 & 0
\end{pmatrix} \in M_{2m}(\C).
\]
Thus, by [HT2, Lemma~3.2],
\[
\begin{pmatrix}
0 & (\lambda\otimes\unit_n - S_n)^{-t}(x\otimes \unit_n)(\lambda\otimes\unit - S_n)^{-1}\\
0 & 0
\end{pmatrix} =
\diff ((\tilde{\lambda}-t\tilde{x})\otimes\unit_n - \tilde{S}_n)^{-1}.
\]
It follows that
\begin{equation*}
\begin{split}
(\id_{m}\otimes&\tr_n)[(\lambda^t\otimes\unit_n - S_n^t)^{-1}(x\otimes \unit_n)(\lambda\otimes\unit - S_n)^{-1}]\\
&=(\Tr_2\otimes\id_m)\Bigg\{
\begin{pmatrix}
0 & (\id_m\otimes\tr_n)[(\lambda\otimes\unit_n - S_n)^{-t}(x\otimes \unit_n)(\lambda\otimes\unit - S_n)^{-1}]\\
0 & 0
\end{pmatrix}
\begin{pmatrix}
0 & 0\\
\unit_m & 0
\end{pmatrix}
\Bigg\}\\
&=(\Tr_2\otimes\id_m)\Bigg\{(\id_{2m}\otimes \tr_n)[(\tilde{\lambda}\otimes\unit_n -\tilde{S_n})^{-1}(x\otimes\unit_n)(\tilde{\lambda}\otimes\unit_n -\tilde{S_n})^{-1}]
\begin{pmatrix}
0 & 0 \\
\unit_m & 0
\end{pmatrix}\Bigg\}\\
&=(\Tr_2\otimes\id_m)\Bigg\{
\diff (\id_{2m}\otimes\tr_n)[((\tilde{\lambda}-t\tilde{x})\otimes\unit_n - \tilde{S}_n)^{-1}]
\begin{pmatrix}
0 & 0\\
\unit_m & 0
\end{pmatrix}
\Bigg\},
\end{split}
\end{equation*}
and
\begin{equation*}
\begin{split}
\E\{(\id_{m}\otimes\tr_n)[(\lambda\otimes\unit_n - S_n)^{-t}&(x\otimes \unit_n)(\lambda\otimes\unit - S_n)^{-1}]\} \\
&=(\Tr_2\otimes\id_m)\Bigg\{
\diff \tilde{G}_n(\tilde{\lambda}-t\tilde{x})
\begin{pmatrix}
0 & 0\\
\unit_m & 0
\end{pmatrix}
\Bigg\}\\
&=- (\Tr_2\otimes\id_m)\Bigg\{
\tilde{G}_n'(\tilde{\lambda})[\tilde{x}]
\begin{pmatrix}
0 & 0\\
\unit_m & 0
\end{pmatrix}
\Bigg\},
\end{split}
\end{equation*}
and the proof is complete. $\endproof$

\vspace{.2cm}

\begin{remark}\label{R}
Inspired by Lemma~\ref{Rn} we define $R(\lambda)$, a 'semicircular analogue' of $R_n(\lambda)$:
\begin{eqnarray*}
R(\lambda) & = &  - \sum_{j=1}^r\sum_{k,l=1}^ma_je_{kl}^{(m)}(\Tr_2\otimes\id_m)\Bigg\{\tilde{G}'(\tilde{\lambda})\Bigg[
\begin{pmatrix}
0 & e_{kl}^{(m)}a_j\\
0 & 0
\end{pmatrix}
\Bigg]
\begin{pmatrix}
0 & 0\\
\unit_m & 0
\end{pmatrix}
\Bigg\}\\
& &  - \sum_{j=1}^r\sum_{k,l=1}^ma_j\cc e_{kl}^{(m)}(\Tr_2\otimes\id_m)\Bigg\{\tilde{G}'(\tilde{\lambda})\Bigg[
\begin{pmatrix}
0 & e_{kl}^{(m)}a_j\cc\\
0 & 0
\end{pmatrix}
\Bigg]
\begin{pmatrix}
0 & 0\\
\unit_m & 0
\end{pmatrix}
\Bigg\}.
\end{eqnarray*}
Then, according to (\ref{diff1}),
\begin{eqnarray*}
\|R_n(\lambda) - R(\lambda)\|& \leq & 2m^2\sum_{j=1}^r\|a_j\|\|(\Tr_2\otimes\id_m)\|\|\tilde{G}_n'(\tilde{\lambda})-\tilde{G}'(\tilde{\lambda})\|\|a_j\| \\
& \leq & \frac{2304\tilde{C}_2 m^2}{n}\Big(\sum_{j=1}^r\|a_j\|^2\Big)(\|\lambda\|+K)^2(\|(\im\lambda)^{-1}\|^8 +\|(\im\lambda)^{-1}\|^6).
\end{eqnarray*}
And again, matrix manipulations reveal that
\begin{equation*}
\begin{split}
R(\lambda) = &  \sum_{j=1}^r\sum_{k,l=1}^m a_j e_{kl}^{(m)}(\id_m\otimes\tau)[(\lambda^t\otimes\unit_\CA - \overline{s})^{-1}(e_{kl}^{(m)}a_j \otimes\unit_\CA)(\lambda\otimes\unit_\CA - s)^{-1}] \\
& +  \sum_{j=1}^r\sum_{k,l=1}^m a_j\cc e_{kl}^{(m)}(\id_m\otimes\tau)[(\lambda^t\otimes\unit_\CA - \overline{s})^{-1}(e_{kl}^{(m)}a_j\cc \otimes\unit_\CA)(\lambda\otimes\unit_\CA - s)^{-1}],
\end{split}
\end{equation*}
where
\[
\overline{s} = \overline{a_0}\otimes \unit_\CA + \sum_{j=1}^r (\overline{a_j}\otimes y_j + \overline{a_j}^*\otimes y_j^*).
\]
In particular,
\begin{equation}\label{normR}
\|R(\lambda)\| \leq 2m^2\Big(\sum_{j=1}^r \|a_j\|^2\Big)\|(\im\lambda)^{-1}\|^2.\end{equation}
The same argument applied to (\ref{MEQa}) gives
\begin{equation}\label{normRn}
\|R_n(\lambda)\| \leq 2m^2\Big(\sum_{j=1}^r \|a_j\|^2\Big)\|(\im\lambda)^{-1}\|^2.
\end{equation}

\end{remark}

\vspace{.2cm}

\begin{remark} Before proceeding any further we mention that all of the
  constants introduced this far, i.e. $C_1$, $C_1'$, $C_2$ and
  $\tilde{C_2}$,  may be expressed in terms of $m, \|a_0\|, \ldots, \|a_{r-1}\|$ and $\|a_r\|$. 
\end{remark}

\vspace{.2cm}

\begin{thm}\label{longthm}
Let $r$ and $m$ be positive integers, let $a_0, a_1, \ldots, a_r$ be
matrices in $M_m(\C)$ with $a_0 = a_0^*$, and for each positive integer $n$
let $Y_1^{(n)}, \ldots , Y_r^{(n)}$ be stochastically independent random
matrices from $\textrm{GRM}^\mathbb{R}(n, \frac{1}{n})$. Furthermore, let $(\mathcal{A}, \tau)$ be a $C^*$-probability space, and let $y_1, \ldots, y_r$ be a circular system in $(\mathcal{A}, \tau)$. Define
\begin{eqnarray*}
  s & = & a_0\otimes \mathbf{1}_{\mathcal{A}} + \sum_{j=1}^r( a_j\otimes y_j + a_j^*\otimes y_j^*),\\
  S_n & = & a_0\otimes \mathbf{1}_n + \sum_{j=1}^r (a_j\otimes Y_j^{(n)}+ a_j^*\otimes (Y_j^{(n)})^*),\\
\end{eqnarray*}
and for $\lambda \in \mathcal{O} = \{\lambda \in M_m(\mathbb{C})| \textrm{Im}\lambda > 0\}$ put
\begin{eqnarray}
  G_n(\lambda) & = &  \mathbb{E}\{(\textrm{id}_m\otimes \textrm{tr}_n)[(\lambda\otimes \mathbf{1}_n-S_n)^{-1}]\}, \label{add1}\\
  G(\lambda) & = & (\textrm{id}_m\otimes \tau)[(\lambda\otimes \mathbf{1}-s)^{-1}], \label{add2}\\
\overline{s} & = & \overline{a_0}\otimes \unit_\CA + \sum_{j=1}^r (\overline{a_j}\otimes y_j + \overline{a_j}^*\otimes y_j^*),\\
R(\lambda) & = &   \sum_{j=1}^r\sum_{k,l=1}^m a_j e_{kl}^{(m)}(\id_m\otimes\tau)[(\lambda^t\otimes\unit_\CA - \overline{s})^{-1}(e_{kl}^{(m)}a_j \otimes\unit_\CA)(\lambda\otimes\unit_\CA - s)^{-1}] \nonumber \\
& &  \label{add3}\\
& + & \sum_{j=1}^r\sum_{k,l=1}^m a_j\cc e_{kl}^{(m)}(\id_m\otimes\tau)[(\lambda^t\otimes\unit_\CA - \overline{s})^{-1}(e_{kl}^{(m)}a_j\cc \otimes\unit_\CA)(\lambda\otimes\unit_\CA - s)^{-1}],\nonumber 
\end{eqnarray}
and
\begin{equation}\label{L}
  L(\lambda) = (\id_m\otimes\tau)[(\lambda\otimes\unit_\CA -s)^{-1}(R(\lambda)G(\lambda)^{-1}\otimes\unit)(\lambda\otimes\unit_\CA -s)^{-1}].
\end{equation} 
Then there is a polynomial $P$ of degree 13 with non-negative coefficients
depending only on $m$,$\|a_0\|, \ldots, \|a_{r-1}\| $ and $\|a_r\|$, such that for any $\lambda \in M_m(\mathbb{C})$ with $\im\lambda > 0$ and for any $n\in \mathbb{N}$
\begin{equation}\label{G-Gn-L}
  \Big\|G(\lambda)-G_n(\lambda)+\frac{1}{n}L(\lambda)\Big\|\leq \frac{1}{n^2}(\|\lambda\| + K')^8 P(\|(\im\lambda)^{-1}\|),
\end{equation}
where $K' = \|a_0\|+ 16\sum_{j=1}^r\|a_j\| +1 = K+1$.
\end{thm}

\proof 
Consider a fixed $n \in \mathbb{N}$, and at first consider an arbitrary $\lambda \in \mathcal{O}_n'$. With $\Lambda_n(\lambda)$ as previously defined,
\begin{eqnarray*}
  \Big\|\Lambda_n(\lambda)-\lambda+\frac{1}{n}R(\lambda)G(\lambda)^{-1}\Big\| &\leq & \Big\|\Lambda_n(\lambda)-\lambda+\frac{1}{n}R_n(\lambda)G_n(\lambda)^{-1}\Big\| +\\
& & \frac{1}{n}\Big\|R_n(\lambda)G_n(\lambda)^{-1}-R(\lambda)G(\lambda)^{-1}\Big\|\\
& \leq & \Big\|\Lambda_n(\lambda)-\lambda+\frac{1}{n}R_n(\lambda)G_n(\lambda)^{-1}\Big\| +\\
& & \frac{1}{n}\Big\|R_n(\lambda)(G_n(\lambda)^{-1}-G(\lambda)^{-1})\Big\| +\\
& & \frac{1}{n}\Big\|(R_n(\lambda )-R(\lambda))G(\lambda)^{-1}\Big\|.
\end{eqnarray*}  
Throughout  the proof we shall make use of the estimates
\begin{equation}
\begin{split}
  \|(\lambda\otimes\unit-s)^{-1}\| \leq  \|(\im \lambda)^{-1}\|,\label{MIEQ6}
\end{split}
\end{equation}
\begin{equation}
\begin{split}
  \|(\Lambda_n(\lambda)\otimes\unit-s)^{-1}\| \leq  2\|(\im \lambda)^{-1}\|,
\end{split}
\end{equation}
\begin{equation}
\begin{split}
  \|G(\lambda)^{-1}\|  \leq  (\|\lambda\| +K)^2\|(\textrm{Im}\lambda)^{-1}\|,\label{G-1}
\end{split}
\end{equation}
\begin{equation}
\begin{split}
 \|G_n(\lambda)^{-1}\|  \leq  (\|\lambda\| +K)^2\|(\textrm{Im}\lambda)^{-1}\|,\label{Gn-2}
\end{split}
\end{equation}
which hold for arbitrary $\lambda \in M_m(\C)$ with $\im\lambda > 0$ (cf. [HT2, Lemma 3.1], (\ref{Lambda_n}), the remarks preceeding [HT2, Lemma~5.4], and (\ref{Gn-1})). 
By Theorem~\ref{MIEQ}
\begin{equation}\label{MIEQ1}
  \Big\|\Lambda_n(\lambda)-\lambda+\frac{1}{n}R_n(\lambda)G_n(\lambda)^{-1}\Big\| \leq \frac{C_1}{n^2}(\|\lambda\|+K)^2\|(\textrm{Im}\lambda)^{-1}\|^5,
\end{equation}
and according to Remark~\ref{R} and (\ref{G-1}) we have
\begin{equation}\label{MIEQ2}
  \frac 1n \Big\|(R_n(\lambda)-R(\lambda))G(\lambda)^{-1}\Big\| \leq
  \frac{2304\,\tilde{C}_2 m^2}{n^2}\Big(\sum_{j=1}^r\|a_j\|^2\Big)(\|\lambda\|+K)^4(\|(\im\lambda)^{-1}\|^9 +\|(\im\lambda)^{-1}\|^7).
\end{equation} 
By (\ref{normRn})
\begin{eqnarray*}
  \|R_n(\lambda)\|  & \leq &  2m^2\Big(\sum_{j=1}^r\|a_j\|^2\Big)\|(\textrm{Im}\lambda)^{-1}\|^2,
\end{eqnarray*}
and from Theorem~\ref{Gn-G}, (\ref{G-1}) and (\ref{Gn-2}) we have
\begin{eqnarray*}
  \|G_n(\lambda)^{-1}-G(\lambda)^{-1}\| & = & \|G_n(\lambda)^{-1}(G(\lambda)-G_n(\lambda))G(\lambda)^{-1}\|\\
& \leq& \|G_n(\lambda)^{-1}\|\|G_n(\lambda)-G(\lambda)\| \|G(\lambda)^{-1}\| \\
& \leq & \frac{C_2}{n}(\|\lambda\|+K)^6(\|(\textrm{Im}\lambda)^{-1}\|^9 +\|(\textrm{Im}\lambda)^{-1}\|^7).\\
\end{eqnarray*}
Therefore
\begin{equation}\label{MIEQ3}
  \frac 1n\Big\|R_n(\lambda)(G_n(\lambda)^{-1}-G(\lambda)^{-1})\Big\|
\leq  \frac{2C_2m^2}{n^2}\Big( \sum_{j=1}^r\|a_j\|^2\Big)(\|\lambda\|+K)^8(\|(\textrm{Im}\lambda)^{-1}\|^{11} +\|(\textrm{Im}\lambda)^{-1}\|^9).
\end{equation}

From (\ref{MIEQ1}), (\ref{MIEQ2}) and (\ref{MIEQ3}) it follows that there
is a polynomial $P_1$ of degree 11 with non-negative coefficients depending
only on  $m$, $\|a_0\|, \ldots, \|a_{r-1}\|$ and $\|a_r\|$, such that for all $\lambda \in \mathcal{O}_n'$
\begin{equation}\label{case1}
  \Big\|\Lambda_n(\lambda)-\lambda+\frac{1}{n}R(\lambda)G(\lambda)^{-1}\Big\| \leq \frac{1}{n^2}(\|\lambda\|+K+1)^8 \cdot 
 P_1(\|(\textrm{Im}\lambda)^{-1}\|).
\end{equation}  
\vspace{.2cm}

We proceed as follows: By (\ref{L})
\[
L(\lambda) = (\id_m\otimes\tau)[(\lambda\otimes\unit_\CA-s)^{-1}(R(\lambda)G(\lambda)^{-1}\otimes\unit_\CA)(\lambda\otimes\unit_\CA-s)^{-1}].
\]
Moreover,
\begin{equation}\label{l-L}
(\lambda\otimes\unit_\CA-s)^{-1} -(\Lambda_n(\lambda)\otimes\unit_\CA-s)^{-1} = 
(\lambda\otimes\unit_\CA-s)^{-1}[(\Lambda_n(\lambda)-\lambda)\otimes\unit_\CA](\Lambda_n(\lambda)\otimes\unit_\CA-s)^{-1}.
\end{equation}
Hence, by Proposition~\ref{biglambda},
\begin{equation}
\begin{split}
 \Big\|G(\lambda)&-G_n(\lambda)+\frac{1}{n}L(\lambda)\Big\|\\
 =&\ \Big\|G(\lambda)-G(\Lambda_n(\lambda))+\frac{1}{n}L(\lambda)\Big\|\\
\leq&\ \Big\|(\lambda\otimes\unit_\CA-s)^{-1} -(\Lambda_n(\lambda)\otimes\unit_\CA-s)^{-1} + \frac{1}{n}(\lambda\otimes\unit_\CA-s)^{-1}(R(\lambda)G(\lambda)^{-1}\otimes\unit_\CA)(\lambda\otimes\unit_\CA-s)^{-1}\Big\|\\
\leq &\ \Big\|(\lambda\otimes\unit_\CA-s)^{-1}[(\Lambda_n(\lambda)-\lambda +\frac 1n R(\lambda)G(\lambda)^{-1})\otimes\unit_\CA](\Lambda_n(\lambda)\otimes\unit_\CA-s)^{-1}\Big\|\\
&+\ \Big\|(\lambda\otimes\unit_\CA-s)^{-1}(\frac 1n R(\lambda)G(\lambda)^{-1}\otimes\unit_\CA)((\lambda\otimes\unit_\CA-s)^{-1} - (\Lambda_n(\lambda)\otimes\unit_\CA-s)^{-1})\Big\|\\
\leq &\ \|(\lambda\otimes\unit_\CA-s)^{-1}\|\|(\Lambda_n(\lambda)\otimes\unit_\CA-s)^{-1}\|\Big\|\Lambda_n(\lambda)-\lambda +\frac 1n R(\lambda)G(\lambda)^{-1}\Big\|\\
&+\ \frac{1}{n}\|(\lambda\otimes\unit_\CA-s)^{-1}\|\|R(\lambda)G(\lambda)^{-1}\|\|(\lambda\otimes\unit_\CA-s)^{-1}-(\Lambda_n(\lambda)\otimes\unit_\CA-s)^{-1}\|\label{MIEQ5}
\end{split}
\end{equation}
From Corollary \ref{CMIEQ} and from (\ref{l-L}) we have
\begin{eqnarray}\label{MIEQ4}
  \|(\lambda\otimes\mathbf{1}-s)^{-1}-(\Lambda_n(\lambda)\otimes\mathbf{1}-s)^{-1}\| & \leq & 2 \|(\textrm{Im}\lambda)^{-1}\|^2\|\lambda - \Lambda_n(\lambda)\|\nonumber\\
& \leq & \frac{2C_1'}{n}(\|\lambda\|+K)^2(\|(\textrm{Im}\lambda)^{-1}\|^7 + \|(\textrm{Im}\lambda)^{-1}\|^5).\nonumber\\
& & 
\end{eqnarray}

By insertion of (\ref{normR}), (\ref{MIEQ6}), (\ref{G-1}) and (\ref{MIEQ4}) into (\ref{MIEQ5}) we find that
\begin{equation*}
\begin{split} 
\Big\|G(\lambda)-G_n(\lambda)&+\frac{1}{n}L(\lambda)\Big\| \leq \frac{2}{n^2}(\|\lambda\|+K+1)^8 P_1( \|(\textrm{Im}\lambda)^{-1}\|)\;\|(\im\lambda)^{-1}\|^2 + \\
 & \frac{4C_1'm^2}{n^2}\Big( \sum_{j=1}^r\|a_j\|^2\Big)(\|\lambda\|+K')^2(\|(\textrm{Im}\lambda)^{-1}\|^7 + \|(\textrm{Im}\lambda)^{-1}\|^5)\|(\textrm{Im}\lambda)^{-1}\|^4.
\end{split}
\end{equation*}
Since $K'=K+1\geq 1$, (\ref{G-Gn-L}) readily follows from this inequality in the case $\lambda\in \CO_n'$.

\vspace{.2cm}

Next, assume that $\lambda \in \mathcal{O}\setminus \mathcal{O}_n'$, i.e.
\[
\frac 12 \leq \frac{C_1'}{n}(\|\lambda\|+K)^2(\|(\textrm{Im}\lambda)^{-1}\|^6 + \|(\textrm{Im}\lambda)^{-1}\|^4).
\]
Then, clearly
\begin{equation*}
\begin{split}
 \Big\|G(\lambda)-G_n&(\lambda)+\frac{1}{n}L(\lambda)\Big\|\\
& \leq   \|G(\lambda)-G_n(\lambda)\| + \frac 1n \|L(\lambda)\|\\
 & \leq    \frac{2C_1'}{n}(\|\lambda\|+K)^2(\|(\textrm{Im}\lambda)^{-1}\|^6 + \|(\textrm{Im}\lambda)^{-1}\|^4)\cdot\Big(\|G(\lambda)-G_n(\lambda)\| + \frac 1n \|L(\lambda)\|\Big),
\end{split}
\end{equation*}
where 
\[
 \|G(\lambda)-G_n(\lambda)\| \leq \frac{C_2}{n}(\|\lambda\|+K)^2(\|(\textrm{Im}\lambda)^{-1}\|^7 + \|(\textrm{Im}\lambda)^{-1}\|^5),
\]
and
\begin{equation}\label{normL}
\|L(\lambda)\| \leq  \|R(\lambda)\|\|G(\lambda)^{-1}\|\|(\im\lambda)^{-1}\|^2 \leq 2m^2\Big( \sum_{j=1}^r\|a_j\|^2\Big)(\|\lambda\|+K)^2\|(\im\lambda)^{-1}\|^5
\end{equation}
Since $K' \geq 1$, it is not hard to see from the estimates above that it is possible to choose a polynomial $P$ having the properties stated in the theorem, such that (\ref{G-Gn-L}) holds in both of the cases $\lambda \in \CO_n'$ and $\lambda \in \CO\setminus\CO_n'$.
$\endproof$

\begin{thm}\label{add}
  For $m\in\N$ put
  \[
  \Om = \{\lambda\in M_m(\C) \,|\,\im \lambda < 0\},
  \]
and for $\lambda\in\Om$ define $G_n(\lambda)$, $G(\lambda)$, $R(\lambda)$ and $L(\lambda)$ by the formulas (\ref{add1}),  (\ref{add2}),  (\ref{add3}) and  (\ref{L}), respectively. Then these four functions taking values in $M_m(\C)$ are well-defined and analytic in $\Om$. Moreover, (\ref{G-Gn-L}) holds for $\lambda\in\Om$ too, with the same polynomial $P$ as in Theorem~\ref{longthm}.
\end{thm}

\proof For $\mu \in \CO$, let $ G_n^{(-)}(\mu)$, $ G^{(-)}(\mu)$, $R^{(-)}(\mu)$ and $ L^{(-)}(\mu)$ be the functions obtained by replacement of $S_n$ and $s$ by $-S_n$ and $-s$ respectively in the formulas (\ref{add1}),  (\ref{add2}),  (\ref{add3}) and  (\ref{L}). Then, for $\lambda\in\Om$ one has that $-\lambda\in\CO$, and 
\begin{eqnarray}
G(\lambda) & = & -G^{(-)}(-\lambda), \label{minus1}\\
G_n(\lambda) & = & -G_n^{(-)}(-\lambda), \label{minus2}\\
R(\lambda) & = & R^{(-)}(-\lambda),\label{minus3}\\
L(\lambda) & = & -L^{(-)}(-\lambda).\label{minus4}
\end{eqnarray}

In particular, $G_n(\lambda)$, $G(\lambda)$, $R(\lambda)$ and $L(\lambda)$ are well-defined and analytic in $\Om$. Moreover, applying (\ref{G-Gn-L}) with $-S_n$ and $-s$ replacing $S_n$ and $s$, respectively, we find that for $\mu\in \CO$
\begin{equation*}
  \Big\| G^{(-)}(\mu)- G_n^{(-)}(\mu)+\frac 1n  L^{(-)}(\mu)  \Big\|\leq \frac{1}{n^2}(\|\mu\| + K')^8 P(\|(\im\mu)^{-1}\|).
\end{equation*}

For $\lambda\in\Om$ we put $\mu = -\lambda \in \CO$ and deduce from (\ref{minus1}),  (\ref{minus2}),  (\ref{minus3}),  (\ref{minus4}) and the estimate above that
\begin{equation*}
  \Big\|G(\lambda)-G_n(\lambda)+\frac{1}{n}L(\lambda)\Big\|\leq \frac{1}{n^2}(\|\lambda\| + K')^8 P(\|(\im\lambda)^{-1}\|). \endproof
\end{equation*}

\section{The spectrum of $S_n$.}

We begin this section with a proof of

\begin{prop}\label{spectra}
Let $(\CA, \tau)$ be a $C^*$-probability space where $\tau$ is a faithful state on $\CA$. Let $r\in \N$, and let $\{y_1, \ldots, y_r\}$ be a circular system in $(\CA, \tau)$. For $m\in \N$ and $a_0, a_1, \ldots, a_r \in M_m(\C)$ with $a_0\cc = a_0$ define
\begin{eqnarray*}
s & = & a_0\otimes\unit_\CA + \sum_{j=0}^r(a_j\otimes y_j + a_j\cc\otimes y_j\cc),\\
\overline{s} & = & \overline{a}_0\otimes\unit_\CA + \sum_{j=0}^r(\overline{a_j}\otimes y_j + \overline{a_j}\cc\otimes y_j\cc),
\end{eqnarray*}
where $\overline{a_j} = (a_j\cc)^t$. Then
\[
\sigma(s) = \sigma(\overline{s}).
\]
\end{prop}

\vspace{.2cm}

\begin{lemma}\label{isomorphisms}
\begin{itemize}
\item[(i)] Let $\{x_1, \ldots, x_s\}$ be a semicircular system w.r.t. a faithful state $\tau$. Then there exists a conjugate linear $\ast$-isomorphism $\pi$ of $C^*(x_1,\ldots, x_s,\unit)$ onto itself, such that
\[
\pi(x_j) = x_j, \;\;\;\;\;\;(1\leq j\leq s),
\]
and $\tau \circ \pi = \overline{\tau}$. 
\item[(ii)] Let $\{y_1, \ldots, y_r\}$ be a circular system w.r.t. a faithful state $\tau$. Then there exists a conjugate linear $\ast$-isomorphism $\rho$ of $C^*(y_1,\ldots, y_r,\unit)$ onto itself, such that
\[
\rho(y_j) = y_j, \;\;\;\;\;\;(1\leq j\leq r),
\]
and $\tau \circ \rho = \overline{\tau}$. 
\end{itemize}
\end{lemma}

\proof
(i) With $\CA = C^*(x_1,\ldots, x_s,\unit)$, let $\CA^c$ be the $C^*$-algebra obtained from $\CA$ in the following way: As a Banach $\ast$-algebra over the reals, $\CA^c$ is identical to $\CA$, but multiplication by complex scalars is changed into multiplication by the complex conjugate scalars. 

For $a\in\CA$ we let $a^c$ denote the corresponding element in $\CA^c$, and we define a faithful state $\tau^c$ on $\CA^c$ by:
\[
\tau^c(a^c) = \overline{\tau(a)}, \;\;\;\;\;\;(a\in\CA). 
\]
As all of the mixed moments of $x_1,\ldots, x_s$ are real numbers, the $\ast$-distribution of $\{x_1^c, \ldots, x_s^c\}$ w.r.t. $\tau^c$ is the same as the $\ast$-distribution of $\{x_1, \ldots, x_s\}$ w.r.t. $\tau$. Furthermore, it is clear that $\unit , x_1^c, \ldots, x_s^c$ generate $\CA^c$. Hence, by [V1, Remark~1.8] there is a unique $\ast$-isomorphism $\phi$ of $\CA$ onto $\CA^c$ such that
\[
\phi(x_j) = x_j^c, \;\;\;\;\;\,(1\leq j\leq s),
\]
and $\tau = \tau^c\circ\phi$. 

Define $\chi : \CA^c\rightarrow \CA$ by
\[
\chi(a^c) = a, \;\;\;\;\;\;(a\in\CA),
\]
and put $\pi = \chi \circ\phi : \CA\rightarrow \CA$. Then $\pi$ has the properties stated in (i). 

(ii) By definition there is a semicircular system $\{x_1, \ldots, x_{2r}\}$ such that
\[
y_j = \frac{x_j+\i x_{j+r}}{\sqrt2}, \;\;\;\;\;\;(1\leq j\leq r).
\]
It is clear that $\unit, x_1, \ldots, x_{2r}$ generate the same $C^*$-algebra $\CA$ as $\unit, y_1, \ldots, y_r$ do. By (i) there is a conjugate linear $\ast$-isomorphism $\pi$ of $\CA$ onto $\CA$ such that
\[
\pi(x_j) = x_j, \;\;\;\;\;\;(1\leq j\leq 2r),
\]
and $\tau\circ\pi = \overline{\tau}$. 

Now, $\{x_1, \ldots, x_r, -x_{r+1}, \ldots, -x_{2r}\}$ is also a semicircular system which, together with $\unit$, generates $\CA$. And then again, by [V1, Remark~1.8] there is a unique $\ast$-automorphism $\psi$ of $\CA$, such that
\begin{eqnarray*}
\psi(x_j) & = & x_j, \;\;\;\;\;\;(1\leq j\leq r),\\
\psi(x_{r+j}) & = & -x_{r+j}, \;\;\;\;\;\;(1\leq j\leq r),
\end{eqnarray*}
and $\tau\circ\psi = \tau$. Finally, let $\rho = \psi\circ\pi$. Then $\rho$ has the properties stated in (ii). $\endproof$

\vspace{.2cm}

\textit{Proof of Proposition~\ref{spectra}}.
Put $\CA_0=C^*(y_1, \ldots, y_r, \unit_\CA)$, and let $\rho :\CA_0 \rightarrow \CA_0$ be the conjugate linear $\ast$-isomorphism provided by Lemma~\ref{isomorphisms}~(ii). $\rho$ extends to a conjugate linear $\ast$-isomorphism 
\[
\rho_m : M_m(\C)\otimes \CA_0 \rightarrow M_m(\C)\otimes \CA_0
\] 
uniquely determined by
\[
\rho_m(x\otimes y) = \overline{x}\otimes \rho(y), \;\;\;\;\;\;(x \in M_m(\C),\; y\in \CA_0).
\]
It is standard to check that for any $v\in M_m(\C)\otimes\CA_0$, $\sigma(\rho_m(v)) = \overline{\sigma(v)}$. In particular, as $s$ is self-adjoint,
\[
\;\;\;\;\;\;\;\;\;\;\;\;\;\sigma(\overline{s}) = \sigma(\rho_m(s)) = \overline{\sigma(s)} = \sigma(s). \;\;\;\;\;\;\;\;\;\;\;\endproof
\]

Let $S_n$, $s$, and $L$ be as defined in Theorem~\ref{longthm}, and for $\lambda \in \C\setminus \R$ define
\begin{eqnarray*}
g_n(\lambda) & = & \E\{(\tr_m\otimes\tr_n)[(\lambda\unit_m\otimes\unit_n - S_n)^{-1}]\},\\
g(\lambda) & = & (\tr_m\otimes\tau)[(\lambda\unit_m\otimes\unit_\CA - s)^{-1}]\},\\
l(\lambda) & = & \tr_m(L(\lambda\unit_m)).
\end{eqnarray*}
By Theorem~\ref{longthm} and Theorem~\ref{add} these are actually well-defined.

\begin{lemma}
Let $\lambda \in \C\setminus \R$ and let $n\in\N$. Then, with $P$ and $K'$ as in  Theorem~\ref{longthm},
\begin{equation}\label{gn-g-l}
\Big|g(\lambda) - g_n(\lambda) +\frac 1n l(\lambda)\Big| \leq \frac{1}{n^2}(|\lambda|+K')^8 P(|\im\lambda|^{-1}).
\end{equation}
\end{lemma}

\proof
This is an easy consequence of Theorem~\ref{longthm} and Theorem~\ref{add}, because
\[
g_n(\lambda) = \tr_m(G_n(\lambda \unit_m)),
\]
and
\[
g(\lambda) = \tr_m(G(\lambda \unit_m)). \endproof
\]

\vspace{.2cm}

Recall from [Ru, Definition~6.7] that a distribution on $\R$ is a linear
functional, $\Lambda : \Ccinf \rightarrow \C$, which is continuous w.r.t. a
certain topology on $\Ccinf$. According to [Ru, Theorem~6.23], for any such
$\Lambda$ there exists a smallest closed set, $F \subseteq \R$, such that
for all $\phi \in \Ccinf$ with $\supp(\phi)\cap F = \emptyset$,
$\Lambda(\phi) = 0$. $F$ is called the support of $\Lambda$ and is denoted
by $\supp(\Lambda)$. We denote by $\Distc$ the set of compactly supported
distributions on $\R$. By [Ru, Theorem~6.24~(d)], every $\Lambda \in
\Distc$ has a natural extension to a linear funtional on $\Cinf$. We let
$\Lambda$ denote this extension as well.

In the proof of Lemma~\ref{distribution} we shall need the following result proved by Tillmann in 1953:

\begin{thm}\label{Tillmann} $[{\rm Ti,  Satz~9~\&~10}]$
  \begin{itemize}
  \item[(i)] Let $\Lambda$ be a distribution on $\R$ with compact support. Define the Stieltjes transform of $\Lambda$, $l:\C\setminus\R \rightarrow \C$, by
    \[
    l(\lambda) = \Lambda\Big(\frac{1}{\lambda -x}\Big), \;\;\;\;\;\;(\lambda \in \C\setminus\R ).
    \]
Then $l$ is analytic in $\C\setminus\R $ and has an analytic continuation to $\C\setminus\supp(\Lambda) $. Moreover,
\begin{itemize}
  \item[(a)] $l(\lambda)\rightarrow 0$, as $|\lambda|\rightarrow \infty$, 
  \item[(b)] there exist a constant $C>0$, an $n\in\N$ and a compact set $K\subset\R$ containing $\supp(\Lambda)$, such that
     \[
     |l(\lambda)| \leq C\cdot\max\{{\rm dist}(\lambda, K)^{-n}, 1\}, \;\;\;\;\;\;(\lambda\in\C\setminus\R),
     \]
   \end{itemize}
and
\begin{itemize}
\item[(c)]  for any $\phi\in\Ccinf$
\[
\Lambda(\phi) = \lim_{y\rightarrow 0^+}\frac{\i}{2\pi}\int_\R \phi(x)[l(x+\i y)-l(x-\i y)]\d x.
\]
    
\end{itemize}

\item[(ii)] Conversely, if $K$ is a compact subset of $\R$, and if $l: \C\setminus K \rightarrow \C$ is an analytic function satisfying (a) and (b) above, then $l$ is the Stieltjes transform of a compactly supported distribution $\Lambda$ on $\R$. Moreover, $\supp(\Lambda)$ is exactly the set of singular points of $l$ in $K$.
\end{itemize}
\end{thm}

\vspace{.2cm}

\begin{lemma}\label{distribution}
There is a distribution $\Lambda \in \Distc$ with $\supp(\Lambda)\subseteq \sigma(s)$, such that for any $\phi \in \Ccinf$
\begin{equation}\label{dist}
\Lambda(\phi) = \lim_{y\rightarrow 0^+}\frac{\i}{2\pi}\int_\R \phi(x)[l(x+\i y)-l(x-\i y)]\d x.
\end{equation}
\end{lemma}

\proof
We show that $l$ satisfies (a) and (b) of Theorem~\ref{Tillmann} (i). By Proposition~\ref{spectra},  $\sigma(\overline{s}) = \sigma(s)$. In addition,
\[
L(\lambda\unit_m) = (\id_m\otimes\tau) [(\lambda\unit_m\otimes\unit_\CA -s)^{-1}(R(\lambda\unit_m) G(\lambda\unit_m)^{-1}\otimes\unit_\CA)(\lambda\unit_m\otimes\unit_\CA -s)^{-1}]
\]
where
\begin{eqnarray*}
R(\lambda\unit_m) & = & \sum_{j=1}^r\sum_{k,l=1}^m a_j e_{kl}^{(m)}(\id_m\otimes\tau)[(\lambda\unit_m\otimes\unit_\CA - \overline{s})^{-1}(e_{kl}^{(m)}a_j \otimes\unit_\CA)(\lambda\unit_m\otimes\unit_\CA - s)^{-1}] \\
& & +  \sum_{j=1}^r\sum_{k,l=1}^m a_j\cc e_{kl}^{(m)}(\id_m\otimes\tau)[(\lambda\unit_m\otimes\unit_\CA - \overline{s})^{-1}(e_{kl}^{(m)}a_j\cc \otimes\unit_\CA)(\lambda\unit_m\otimes\unit_\CA - s)^{-1}]
\end{eqnarray*}
and
\[
G(\lambda\unit_m)^{-1} = \lambda\unit_m -a_0 - \sum_{j=1}^r(a_jG(\lambda\unit_m)a_j\cc + a_j\cc G(\lambda\unit_m)a_j).
\]
It follows that $L$ and $l$ have analytic continuations to $\C\setminus\sigma(s)$.

Moreover,
\[
|l(\lambda)| \leq \|L(\lambda\unit_m)\|\leq \|R(\lambda\unit_m)\|\|G(\lambda\unit_m)^{-1}\|\|(\lambda\unit_m\otimes\unit_\CA -s)^{-1}\|^2,
\]
where, for $|\lambda|>\|s\|$, 
\begin{eqnarray*}
\|R(\lambda\unit_m)\| & \leq & 2m^2\Big(\sum_{j=1}^r\|a_j\|^2\Big)\|(\lambda\unit_m\otimes\unit_\CA -s)^{-1}\|^2\\
& \leq &  2m^2\Big(\sum_{j=1}^r\|a_j\|^2\Big)\Big( \frac{1}{|\lambda| -\|s\|} \Big)^2,
\end{eqnarray*}
and
\begin{eqnarray*}
\|G(\lambda\unit_m)^{-1}\| & = & \Big\|\lambda\unit_m -a_0 - \sum_{j=1}^r(a_jG(\lambda\unit_m)a_j\cc + a_j\cc G(\lambda\unit_m)a_j)\Big\|\\
& \leq & |\lambda| + \|a_0\| +2\Big(\sum_{j=1}^r\|a_j\|^2\Big) \frac{1}{|\lambda| -\|s\|}.
\end{eqnarray*}
It follows that
\begin{equation}\label{absl}
|l(\lambda)|\leq O\Big(\frac{1}{|\lambda|^3}\Big), \;\; {\rm as} \;\; |\lambda| \rightarrow \infty. 
\end{equation}
Hence, (a) is satisfied.

By (\ref{normL}), for $\lambda\in \C\setminus\R$
\begin{equation}\label{final}
|l(\lambda)| \leq \|L(\lambda\unit_m)\| \leq 2m^2 \Big(\sum_{j=1}^r\|a_j\|^2\Big)(|\lambda| + K')^2|(\im\lambda)^{-1}|^5.
\end{equation}

Now, choose $a, \, b\in\R$, $a < b$, such that $\sigma(s) \subseteq [a,b]$,
put $K=[a-1,b+1]$, and put
\[
D = \{\lambda\in\C\, |\, 0 < {\rm dist}(\lambda, K)\leq 1 \}.
\]
Note that for $\lambda \in D$, either ${\rm dist}(\lambda, K) = |\im
\lambda|$ or ${\rm dist}(\lambda, \sigma(s)) \geq 1$. Hence, by (\ref{final}) and the fact that $l$ is bounded on compact
subsets of $\C\setminus \sigma(s)$, we find that for some constant $C<\infty$:
\[
\forall \lambda \in D:\;\; |l(\lambda)| \leq C\cdot \max\{ {\rm
  dist}(\lambda, K)^{-5}, 1\} = C\cdot {\rm
  dist}(\lambda, K)^{-5} . 
\]

By (\ref{absl}) $l$ is bounded on $\C\setminus D$. Hence, if $C$
is chosen sufficiently large, then (b) holds with $K=[a-1,b+1]$ and $n=5$, and the lemma follows from Theorem~\ref{Tillmann}. $\endproof$

\vspace{.2cm}

\begin{thm}\label{distthm}
Let $S_n$, $s$, $P$ and $K'$ be as in Theorem~\ref{longthm}, and let $\Lambda \in \Distc$ be as in Lemma~\ref{distribution}. Then for any $\phi \in \Ccinf$:
\begin{equation}
\begin{split}
\Big|\E\{(\tr_m\otimes\tr_n)&\phi(S_n)\}-(\tr_m\otimes\tau)\phi(s) -\frac 1n \Lambda(\phi)\Big|\\
&\leq\frac{1}{n^2}\Big(\frac{2^7}{13!\pi}\int_{\R\times\R}|((1+D)^{14}\phi)(x)|(|x|+\sqrt{2}t+K')^8 Q(t)\e^{-t} \d m(x)\otimes\d m(t)\Big),
\end{split}
\end{equation}
where $D=\frac{\d }{\d x}$, and $Q$ is the polynomial of degree less than 13 defined by $Q(t) = P(t^{-1})t^{13}$. In particular,
\begin{equation}
\E\{(\tr_m\otimes\tr_n)\phi(S_n)\}=(\tr_m\otimes\tau)\phi(s) +\frac 1n \Lambda(\phi) + O\Big(\frac{1}{n^2}\Big).
\end{equation}
\end{thm}

\proof
By the Riesz representation theorem there are unique probability measures $\mu_n$ and $\mu$ on $(\R, \B)$ such that for any $\psi \in C_0(\R)$
\[
\int_\R \psi \;\d \mu_n = \E\{(\tr_m\otimes\tr_n)\psi(S_n)\},
\]
and 
\[
\int_\R \psi \;\d \mu = (\tr_m\otimes\tau)\psi(s).
\]
In particular, for $\lambda \in \C\setminus\R$ we have
\[
g_n(\lambda) = \int_\R\frac{1}{\lambda -x}\d \mu_n(x),
\]
and
\[
g(\lambda) = \int_\R\frac{1}{\lambda -x}\d \mu(x).
\]
By the inverse Stieltjes transform
\[
\d \mu_n(x) = \lim_{y\rightarrow 0^+} \Big(-\frac 1\pi \im(g_n(x+\i y))\d x\Big) = \lim_{y\rightarrow 0^+} \Big(\frac{\i}{2\pi}(g_n(x+\i y)-g_n(x-\i y))\d x\Big),
\]
and
\[
\d \mu(x) = \lim_{y\rightarrow 0^+} \Big(-\frac 1\pi \im(g(x+\i y))\d x\Big) = \lim_{y\rightarrow 0^+} \Big(\frac{\i}{2\pi}(g(x+\i y)-g(x-\i y))\d x\Big),
\]
in the sense of weak convergence of probability measures on $(\R, \B)$. Combining these observations with Lemma~\ref{distribution} we find that for any $\phi \in \Ccinf$
\begin{equation*}
\begin{split}
\E\{(\tr_m\otimes\tr_n)\phi(S_n)\}-(\tr_m\otimes\tau)\phi(s) -\frac 1n \Lambda(\phi) =\\
\lim_{y\rightarrow 0^+}\frac{\i}{2\pi}\int_\R \phi(x) [g_n(x+\i y)-g(x+\i y)& - \frac 1n l(x+ \i y) -\\
&(g_n(x-\i y)-g(x-\i y) - \frac 1n l(x-\i y))]\d x.
\end{split}
\end{equation*}
Thus, with
\[
r_n(\lambda) = g_n(\lambda)-g(\lambda) - \frac 1n l(\lambda), \;\;\;\;\;\;(\lambda \in \C\setminus\R),
\]
\begin{equation}\label{spectrum1}
\begin{split}
\Big|\E\{(\tr_m\otimes\tr_n)&\phi(S_n)\}-(\tr_m\otimes\tau)\phi(s) -\frac 1n \Lambda(\phi)\Big| \\
&\leq\frac{1}{2\pi}\limsup_{y\rightarrow 0^+}\Big(\Big|\int_\R \phi(x) r_n(x+\i y)\d x\Big| + \Big|\int_\R \phi(x) r_n(x-\i y)\d x\Big|\Big).
\end{split}
\end{equation}
The rest of the proof follows the lines of [HT2, Proof of Theorem~6.2]. For $\im\lambda \neq 0$ and $p\in \N$ define
\[
I_p(\lambda) = \frac{1}{(p-1)!}\int_\R r_n(\lambda +t) t^{p-1}\e^{-t}\d t.
\]
By (\ref{gn-g-l}), if $\eps > 0$ and $|\im\lambda| >\eps$, then for every $t >0$\[
|r_n(\lambda +t)|\leq \frac{1}{n^2}(|\lambda +t|+K')^6 P(|\im(\lambda +t)|^{-1}) \leq \frac{1}{n^2}(|\lambda| +t + K')^8 P(\eps^{-1}).
\]

Therefore $I_p(\lambda)$ is well-defined, and it is standard to check that $I_p$ is analytic. 

\vspace{.2cm}

Now, consider a fixed $\phi \in \Ccinf$, and let $y>0$. After $p$ steps of integration by parts we obtain:
\[
\int_\R \phi(x)r_n(x\pm \i y)\d x = \int_\R ((1+D)^p\phi)(x)I_p(x\pm \i y)\d x,
\]
where $D= \frac{\d}{\d x}$. 

We are going to estimate $|I_{14}(\lambda)|$. Let's start with the case $\im\lambda >0$ and define
\[
F(z) = \frac{1}{13!}r_n(\lambda+z)z^{13}\e^{-z}, \;\;\;\;\;\;(\im z > - \im\lambda).
\]
$F$ is analytic, so for any $r > 0$
\[
\int_{[0, r]}F(z)\d z + \int_{[r, r+\i r]}F(z)\d z  + \int_{[r+\i r, 0]}F(z)\d z =0.
\]
The second term in the expression above tends to zero, as $r$ goes to infinity. Indeed,
\begin{equation*}
\begin{split}
\Big|\int_{[r, r+\i r]}F(z)&\d z \Big|  = \Big|\int_0^r F(r+\i t)\d t \Big|\\
\leq & \;\frac{2^{13}}{13!} \int_0^r|r_n(\lambda +r+\i t)|r^{13}\e^{-r}\d t\\
 \leq &\; \frac{2^{13}}{13!} \int_0^r\Big(|g_n(\lambda +r+\i t)|+ |g(\lambda +r+\i t)| +\frac 1n |l(\lambda +r+\i t)|\Big)r^{13}\e^{-r}\d t\\
\leq & \;\frac{2^{13}}{13!} \int_0^r\Bigg(\frac{2}{|\im\lambda +t|}+\frac 1n 2m^2(\sum_{j=1}^r\|a_j\|^2)(|\lambda +r+\i t|+K)^2\frac{1}{|\im\lambda + t|^5}\Bigg)r^{13}\e^{-r} \d t\\
\leq & \;\frac{2^{13}}{13!} \Bigg(\frac{2}{|\im\lambda|}+\frac 1n 2m^2(\sum_{j=1}^r\|a_j\|^2)(|\lambda| +2r+K)^2\frac{1}{|\im\lambda|^5}\Bigg)r^{14}\e^{-r}\\
\rightarrow & \;0, \;\; {\rm as} \;\;r\rightarrow \infty.
\end{split}
\end{equation*}
Consequently,
\begin{eqnarray*}
I_{14}(\lambda) &=& \lim_{r \rightarrow \infty}\int_{[0, r+\i r]}F(z)\d z\\
& = & \int_0^\infty F((1+\i)t)(1+\i)\d t\\
& = & \frac{1}{13!}  \int_0^\infty r_n(\lambda + (1+\i)t)(1+\i)^{14}t^{13}\e^{-(1+\i)t}\d t,
\end{eqnarray*}
and by (\ref{gn-g-l})
\begin{eqnarray*}
|I_{14}(\lambda)| & \leq & \frac{2^7}{13!n^2}\int_0^\infty (|\lambda + (1+\i )t|+K')^8 P(|\im\lambda + t|^{-1}) t^{13}\e^{-t}\d t\\
& \leq & \frac{2^7}{13!n^2}\int_0^\infty (|\lambda|+ \sqrt{2}t +K')^8 P(t^{-1})t^{13}\e^{-t} \d t\\
& = &  \frac{2^7}{13!n^2}\int_0^\infty (|\lambda|+ \sqrt{2}t +K')^8 Q(t)\e^{-t} \d t
\end{eqnarray*}
with $Q(t) =  P(t^{-1})t^{13}$. The case $\im\lambda <0$ is treated similarly, and we obtain the same upper bound as the one above. Thus, for every $y>0$
\begin{equation*}
\begin{split}
 \Big|\int_\R \phi(x)&r_n(x\pm \i y) \d x\Big|\\
&\leq \;\frac{1}{n^2} \Big(\frac{2^7}{13!}\int_{\R\times\R}|((1+D)^{14}\phi)(x)|(|x\pm \i y|+\sqrt{2}t+K')^8 Q(t)\e^{-t} \d m(x)\otimes\d m(t)\Big),
\end{split}
\end{equation*}
and by dominated convergence,  
\begin{equation*}
\begin{split}
\limsup_{y\rightarrow 0^+} \Big|\int_\R &\phi(x)r_n(x\pm \i y) \d x\Big| \\
&\leq \frac{1}{n^2} \Big(\frac{2^7}{13!}\int_{\R\times\R}|((1+D)^{14}\phi)(x)|(|x|+\sqrt{2}t+K')^8 Q(t)\e^{-t} \d m(x)\otimes\d m(t)\Big).
\end{split}
\end{equation*}
By (\ref{spectrum1}) this completes the proof. $\endproof$

\vspace{.2cm}

An application of Remark~\ref{psi} and Lemma~\ref{sumnorm} yields (as in
[HT2, Proof of Proposition~4.7]) the following estimate:

\begin{prop}\label{var}
Let $\psi : \R\rightarrow \C$ be a $C^1$-function with compact support. Then
\[
\V\{(\tr_m\otimes\tr_n)\psi(S_n)\} \leq \frac{4}{n^2} \Big(\sum_{j=1}^r \|a_j\|^2\Big)\, \E\{(\tr_m\otimes\tr_n)|\psi'|^2(S_n)\}.
\]
\end{prop} 

\vspace{.2cm}

Finally, as in [HT2, Proof of Lemma~6.3] and [HT2, Proof of  Theorem~6.4], we may combine Theorem~\ref{distthm} and Proposition~\ref{var} to get:

\begin{thm}\label{spectrum}
Let $m\in\N$, let $a_0, \ldots, a_r \in M_m(\C)$ with $a_0\cc = a_0$, and let $S_n$ and $s$ be as defined in Theorem~\ref{longthm}. Then for any $\eps > 0$ and for almost all $\omega \in \Omega$,
\[
\sigma(S_n(\omega)) \subseteq \sigma(s) + (-\eps, \eps),
\]
eventually as $n\rightarrow \infty$.
\end{thm}

\vspace{.2cm}

\begin{remark}
Let $r$ and $s$ be non-negative integers with $r+s \geq 1$, and for each
$n\in \N$, let  $X_1^{(n)}, \ldots, X_{r+s}^{(n)}$ be stochastically
independent random matrices such that  $ X_1^{(n)}, \ldots,  X_r^{(n)}\in
\GOE(n, \frac 1n)$ and $ X_{r+1}^{(n)},\ldots, X_{r+s}^{(n)} \in \GOES(n,
\frac 1n)$. Furthermore, let $a_0, b_1, \ldots, b_{r+s} \in M_m(\C)_{sa}$, and put
\[
S_n = a_0\otimes\unit_n + \sum_{j=1}^{r+s} b_j\otimes X_j^{(n)}.
\]
For $\{x_1, \ldots, x_{r+s}\}$ a semicircular system in $(\CA, \tau)$ put
\[
s = a_0\otimes\unit_\CA + \sum_{j=1}^{r+s} b_j\otimes x_j.
\]
As mentioned in the introduction, the $X_j^{(n)}$'s may be expressed in terms of $r+s$ stochastically independent random matrices $Y_1^{(n)}, \ldots , Y_{r+s}^{(n)}\in \GRMR(n, \frac 1n )$, and for suitable $a_1, \ldots, a_r \in M_m(\C)$
\[
 S_n  =  a_0\otimes \mathbf{1}_n + \sum_{j=1}^r (a_j\otimes Y_j^{(n)}+ a_j^*\otimes {Y_j^{(n)}}^*).
\]
We may also assume that there is a circular system $\{y_1, \ldots, y_{r+s}\}$ in $(\CA, \tau)$ such that
\begin{eqnarray*}
x_j =  \frac{1}{\sqrt{2}}(y_j +y_j\cc), &\;\;& (1\leq j \leq r),\\
x_j =  \frac{1}{\i\sqrt{2}}(y_j -y_j\cc) , &\;\;& (r+1\leq j \leq r+s).
\end{eqnarray*}
Then, with $a_1, \ldots, a_r$ as above,
\[
s = a_0\otimes \mathbf{1}_{\mathcal{A}} + \sum_{j=1}^{r+s}( a_j\otimes y_j + a_j^*\otimes y_j^*),
\]
and it follows from Theorem~\ref{spectrum} that for any $\eps > 0$
\[
\sigma(S_n(\omega))\subseteq \sigma(s) + (-\eps, \eps),
\]
eventually as $n\rightarrow \infty$ for almost all $\omega \in \Omega$.
Thus, by the proof of [HT2, Proposition~7.3] we have:
\end{remark}

\begin{prop}
Let $r$ and $s$ be non-negative integers with $r+s \geq 1$, and for each $n\in \N$, let  $X_1^{(n)}, \ldots, X_{r+s}^{(n)}$ be stochastically independent random matrices such that  $ X_1^{(n)}, \ldots,  X_r^{(n)}\in \GOE(n, \frac 1n)$ and $ X_{r+1}^{(n)},\ldots, X_{r+s}^{(n)} \in \GOES(n, \frac 1n)$. 
Furthermore, let $\{x_1, \ldots, x_{r+s}\}$ be a semicircular system in a $C^*$-probability space $(\CA, \tau)$ with $\tau$ a faithful state on $\CA$. Then there is a $P$-null set $N' \subseteq \Omega$ such that for any non-commutative polynomial $p$ in $r+s$ variables and for every $\omega \in \Omega \setminus N'$:
\begin{equation}
  \limsup_{n \rightarrow \infty}\|p(X_1^{(n)}(\omega), \ldots , X_{r+s}^{(n)}(\omega))\| \leq \|p(x_1, \ldots , x_{r+s})\|.
\end{equation}
\end{prop}

\section{Almost sure convergence of mixed moments.}

Throughout this section let $r$ and $s$ be fixed numbers in $\N_0$ with
$r+s\geq 1$, and for each $n\in\N$, let $X_1^{(n)}, \ldots, X_{r+s}^{(n)}$ be stochastically independent random matrices such that  $ X_1^{(n)}, \ldots,  X_r^{(n)}\in \GOE(n, \frac 1n)$ and $ X_{r+1}^{(n)},\ldots, X_{r+s}^{(n)} \in \GOES(n, \frac 1n)$.
  
Also, let $(\CA,\tau)$ denote a $C\cc$-probability space, and let $x_1, \ldots, x_{r+s}$ be a semicircular system in $(\CA,\tau)$.

\vspace{.2cm}

The aim of this section is to prove:

\begin{prop}
There is a $P$-null set $N'' \subseteq \Omega$ such that for any non-commutative polynomial $p$ in $r+s$ variables and for every $\omega \in \Omega \setminus N''$:
\begin{equation}
  \liminf_{n \rightarrow \infty}\|p(X_1^{(n)}(\omega), \ldots , X_{r+s}^{(n)}(\omega))\| \geq \|p(x_1, \ldots , x_{r+s})\|.
\end{equation}
\end{prop}

\vspace{.2cm}

From [HT2, Proof of Lemma~7.2] it is clear that this proposition follows from the following theorem:

\begin{thm}\label{lim} 
For every $p \in \mathbb{C}\< X_1, \ldots, X_{r+s} \>$ 
\begin{equation}
  \lim_{n \rightarrow \infty} \tr_n(p(X_1^{(n)}(\omega), \ldots , X_{r+s}^{(n)}(\omega))) = \tau(p(x_1, \ldots, x_{r+s}))
\end{equation}
holds for almost all $\omega \in \Omega$.
\end{thm}

\vspace{.2cm}
Note that Theorem~\ref{lim} is the strong version of the following result due to Voiculescu:

\begin{thm}$[$V2, Theorem~2.3$]$\label{Voiculescu}
Let $r$, $s$ and $X_1^{(n)}, \ldots, X_{r+s}^{(n)}$ be as defined above. Then 
\begin{itemize}
  \item[(a)] the sets $\{X_j^{(n)}\}_{1\leq j\leq r+s}$ are asymptotically free as $n \rightarrow \infty$, and
    \item[(b)] for all $j \in \{1, \ldots , r+s\}$ and for all $k \in \N$ 
\begin{equation}
  \lim_{n\rightarrow \infty}\E\{\tr_n((X_j^{(n)})^k)\} = \frac{1}{2\pi}\int_{-2}^{2} x^k\sqrt{4-x^2}dx.
\end{equation}
\end{itemize}

In particular, for any $p \in \mathbb{C}\langle X_1, \ldots, X_{r+s} \rangle$
\begin{equation}
  \lim_{n\rightarrow \infty}\E\{\tr_n(p(X_1^{(n)},\ldots ,X_{r+s}^{(n)}))\} =\tau(p(x_1, \ldots, x_{r+s})).
\end{equation}
\end{thm}

\vspace{.2cm}

Actually [V2, Theorem~2.3] is slightly different from
Theorem~\ref{Voiculescu} above, since in [V2] the diagonal elements of the
first $r$ random matrices $X_1^{(n)}, \ldots, X_{r}^{(n)}$ have variance
$\frac 1n$ and not $\frac 2n$ as in the
definition of $\GOE(n,\frac 1n)$. However, when one works out the details
of the proof sketched in [V2], it is easily seen that this change of the
diagonal entries does not affect the validity of the result.

In order to prove that Theorem~\ref{Voiculescu} implies Theorem~\ref{lim}
we shall make use of the following two lemmas:

\begin{lemma}\label{norm}
Let $k\in \mathbb{N}$. Then there is a constant $C_1(k)<\infty$ such that for all $n \in \mathbb{N}$ and for all $X^{(n)} \in \textrm{SGRM}(n, \frac{1}{n})$
\begin{equation}
  \E\{\|X^{(n)}\|^k\}\leq C_1(k).
\end{equation}
\end{lemma}

\proof According to [HT2, Lemma 5.1] the constant 4 dominates
$\E\{\|X^{(n)}\|\}$ for all $n\in \N$. Hence, we shall concentrate on the case $k\geq 2$. Let $n \in \N$, and let $X^{(n)} \in \SGRM(n, \frac{1}{n})$. Let $\lmax(X^{(n)})$ (respectively $\lmin(X^{(n)})$) denote the largest (respectively the smallest) eigenvalue of $X^{(n)}$. Then for all $\eps > 0$ we have (cf. [HT1, Proof of Lemma~3.3]) 
\begin{eqnarray*}
  P(\lmax(X^{(n)}) > 2+\eps ) & \leq & n\cdot \exp{\Big(-\frac{n\eps^2}{2}\Big)},\\
  P(\lmin(X^{(n)}) < -(2+\eps) ) & \leq & n\cdot \exp{\Big(-\frac{n\eps^2}{2}\Big)},
\end{eqnarray*}
where the last equality follows from the first and the fact that $-X^{(n)} \in \SGRM(n, \frac{1}{n})$. Hence,
\begin{equation}\label{est1}
  P(\|X^{(n)}\| > 2+\eps) \leq 2n\cdot \exp{\Big(-\frac{n\eps^2}{2}\Big)}.
\end{equation}
Put 
\begin{equation}
  \eps_0 = \sqrt\frac{2\log{2n}}{n},
\end{equation}
and note that
\begin{equation}\label{id1}
  2n\cdot \exp{\Big(-\frac{n\eps_0^2}{2}\Big)} = 1.
\end{equation}
Define $F: \R\rightarrow [0, 1]$ by
\begin{equation}
  F(t) = P(\|X^{(n)}\|\leq t), \;\;\; (t\in\R).
\end{equation}

Then from the estimate (\ref{est1}) it follows that for all $\eps >0$
\begin{equation}
  F(2+\eps) \geq 1 -  2n\cdot\exp{\Big(-\frac{n\eps^2}{2}\Big)},
\end{equation}
and carrying out integration by parts (cf. [Fe, Lemma~V.6.1]) we obtain
\begin{eqnarray*}
  \E\{\|X^{(n)}\|^k\} & = & \int_0^\infty t^k \d F(t)\\
  & = & k\int_0^\infty t^{k-1}(1-F(t))\d t\\
  & = & k\int_0^{2+\eps_0} t^{k-1}(1-F(t))\d t + k\int_{2+\eps_0}^\infty t^{k-1}(1-F(t))\d t \\
  & \leq & (2+\eps_0)^k + k \int_{\eps_0}^\infty (2+t)^{k-1} 2n \exp{\Big(-\frac{nt^2}{2}\Big)}\d t.
\end{eqnarray*}
Now, using (\ref{id1}) we get
\begin{equation*}
\begin{split}
k \int_{\eps_0}^\infty (2+t)^{k-1}2n\exp&{\Big(-\frac{nt^2}{2}\Big)}\d t \\
&= k \int_{\eps_0}^\infty (2+t)^{k-1}2n\exp{\Big(-\frac{n(t-\eps_0)^2}{2}-\frac{n\eps_0^2}{2}-n\eps_0(t-\eps_0)\Big)}\d t \\
&\leq  k \; 2n\exp{\Big(-\frac{n\eps_0^2}{2}\Big)} \int_{\eps_0}^\infty (2+t)^{k-1}\exp{\Big(-\frac{n(t-\eps_0)^2}{2}\Big)}\d t\\
&=  k \int_0^\infty (2+t+\eps_0)^{k-1}\exp{\Big(-\frac{nt^2}{2}\Big)}\d t \\
&\leq  k \int_0^\infty (2+t+\eps_0)^{k-1}\exp{\Big(-\frac{t^2}{2}\Big)}\d t.  
\end{split}
\end{equation*}
It is easily shown that the function $g: [\frac{1}{2}, \infty[\rightarrow \mathbb{R}$ defined by
\begin{equation}
  g(x) = \sqrt{\frac{2 \log{2x}}{x}}, \;\;\;\;(x\geq\frac 12)
\end{equation}
attains its maximum at $x=\e$, and this maximum is $2\e^{-\frac 12}$ which is less than 2. Hence, 
\[
  \E\{\|X^{(n)}\|^k\} \leq  4^k +  k \int_0^\infty (4+t)^{k-1}\exp{\Big(-\frac{t^2}{2}\Big)}\d t,
\]
and the lemma follows. $\endproof$

\vspace{.2cm}
 
\begin{lemma}\label{variance}
Let  $d \in \N$, let $i_1, \ldots , i_d \in \{1, \ldots , r+s\}$, and let $n \in \N$. Define $f: M_n(\C)^{r+s}\rightarrow \C$ by
\begin{equation}
  f(v_1, \ldots , v_{r+s}) = \tr_n(v_{i_1}\cdot \cdots \cdot v_{i_d}), \;\; (v_1, \ldots , v_{r+s} \in M_n(\mathbb{C})).
\end{equation}
Then there is a constant $C_2(d) > 0$ (independent of $n$) such that
\begin{equation}
  \V\{f(X_1^{(n)},\ldots ,X_{r+s}^{(n)})\}\leq \frac{C_2(d)}{n^2}.
\end{equation}
\end{lemma}

\proof From arguments similar to those presented in Remark~\ref{psi} it follows that
\begin{equation}\label{Poincare}
  \V\{f(X_1^{(n)},\ldots ,X_{r+s}^{(n)})\} \leq \frac{1}{2n}\E\{\|(\grad f)(X_1^{(n)}, \ldots,X_{r+s}^{(n)} )\|_e^2\}.
\end{equation}
Now, let $v = (v_1, \ldots, v_{r+s}) \in M_n(\C)^{r+s}$, and let $w = (w_1, \ldots, w_{r+s}) \in M_n(\C)^{r+s}$ with $\|w\|_e=1$. By the Cauchy-Schwartz inequality:
\begin{equation*}
\begin{split}
\Big|\diff &f(v+tw)\Big|  =  \frac{1}{n}\Big|\diff\Tr_n((v_{i_1}+tw_{i_1})\cdot \cdots \cdot (v_{i_d}+tw_{i_d}))\Big|\\
& \leq  \frac{1}{n}\|\unit_n\|_{2, \Tr_n}(\|w_{i_1}v_{i_2}\cdots v_{i_d}\|_{2, \Tr_n}+ \|v_{i_1}w_{i_2}\cdots v_{i_d}\|_{2, \Tr_n}+\cdots + \|v_{i_1}v_{i_2}\cdots w_{i_d}\|_{2, \Tr_n})\\
& \leq  \frac{1}{\sqrt{n}}(\|w_{i_1}\|_{2, \Tr_n}\|v_{i_2}\cdots v_{i_d}\|+\|v_{i_1}\|\|w_{i_2}\|_{2, \Tr_n}\|v_{i_3}\cdots v_{i_d}\|+ \cdots + \|v_{i_1}\cdots v_{i_{d-1}}\|\|w_{i_d}\|_{2, \Tr_n})
\end{split}
\end{equation*}
Hence, with $M = \max_{1\leq j \leq r+s}\|v_j\|$,
\begin{eqnarray*}
\Big|\diff f(v+tw)\Big| & \leq &  \frac{1}{\sqrt{n}}M^{d-1}\sum_{j=1}^d\|w_{i_j}\|_{2, \Tr_n}\\
& \leq &  \frac{1}{\sqrt{n}}M^{d-1}d,
\end{eqnarray*}
and since $w$ was arbitrary,
\begin{equation}
  \|\grad f\|_e^2\leq \frac{1}{n}M^{2(d-1)}d^2.
\end{equation}
It follows from (\ref{Poincare}) that
\begin{eqnarray*}
  \V\{f(X_1^{(n)},\ldots ,X_{r+s}^{(n)})\}& \leq & \frac{d^2}{2n^2}\E\{(\max_{1\leq j \leq r+s}\|X_j^{(n)}\|)^{2(d-1)}\}\\
  & \leq & \frac{d^2}{2n^2}\sum_{j=1}^{ r+s}\E\{\|X_j^{(n)}\|^{2(d-1)}\}.
\end{eqnarray*}
For any $X\in \GOE(n,\frac{1}{n})$ we may choose $Y \in \SGRM(n,\frac{1}{n})$ such that 
\begin{equation}
  X=\frac{1}{\sqrt{2}}(Y+\overline{Y}),
\end{equation}
Similarly, when $X\in \GOES(n,\frac{1}{n})$ we may choose $Y \in \SGRM(n,\frac{1}{n})$ such that 
\begin{equation}
  X=\frac{1}{\sqrt{2}}(Y-\overline{Y}).
\end{equation}
Then in both cases we have
\begin{equation}
  \E\{\|X\|^{2(d-1)}\} \leq 2^{d-1}\E\{\|Y\|^{2(d-1)}\},
\end{equation} 
and applying Lemma~\ref{norm} we obtain the desired estimate with 
\[
  C_2(d) = d^2 (r+s) 2^{d-2}\cdot C_1(2d-2).\;\;\;\;\;\;\endproof
\] 

\vspace{.2cm}

\textit{Proof of Theorem~\ref{lim}.} Let $d \in \N$, and let $i_1, \ldots , i_d \in \{1, \ldots , r+s\}$. For each $n \in \N$ define a complex random variable $Z_n$ by
\begin{equation}
  Z_n = \tr_n(X_{i_1}^{(n)}\cdot \cdots \cdot  X_{i_d}^{(n)})- \E\{\tr_n(X_{i_1}^{(n)}\cdot \cdots \cdot  X_{i_d}^{(n)})\}.
\end{equation}
By the Borel-Cantelli Lemma, if
\begin{equation}\label{sum}
  \sum_{n=1}^\infty P(|Z_n|> n^{-1/3}) < \infty,
\end{equation}
then 
\begin{equation}
  P(|Z_n|\leq n^{-1/3}, \textrm{eventually as } n \rightarrow \infty) = 1.
\end{equation}
In particular, the sequence $(Z_n)_{n=1}^\infty$ tends to zero almost surely, and it follows from Theorem~\ref{Voiculescu} that 
\begin{equation}\label{last}
\lim_{n\rightarrow\infty}\tr_n(X_{i_1}^{(n)}\cdots X_{i_d}^{(n)}) = \tau(x_{i_1}\cdots x_{i_d}),
\end{equation}
almost surely. Consequently, Theorem~\ref{lim} holds. To prove (\ref{sum}), apply Chebychev's inequality and Lemma~\ref{variance} as follows:
\begin{eqnarray*}
   P(|Z_n|> n^{-1/3}) & \leq & n^{2/3} \E\{|Z_n|^2\}\\
   & \leq & n^{2/3} \frac{C_2(d)}{n^2}\\
   & = & n^{-4/3}C_2(d).\;\;\;\;\;\;\endproof
\end{eqnarray*}

\vspace{.2cm}

\begin{remark} The above proof of Theorem~\ref{lim} follows the main lines
   of the proof of the corresponding result for the $\SGRM(n,\frac
   1n)$-case sketched in [Pi, Proof of Theorem~9.9.3].
\end{remark}

\section{The symplectic case.}

We shall use the results of the previous sections to prove:

\begin{thm}\label{symplectic}
Let $r$, $s \in \N_0$ with $r+s\geq 1$ , and for each $n\in \N$, let
$X_1^{(n)}, \ldots, X_{r+s}^{(n)}$ be stochastically independent random
matrices defined on $(\Omega, \CF, P)$ such that $ X_1^{(n)}, \ldots,
X_r^{(n)}\in \GSE(n, \frac 1n)$ and $ X_{r+1}^{(n)},\ldots, X_{r+s}^{(n)}
\in \GSES(n, \frac 1n)$. Furthermore, let $(\CA, \tau)$ be a
$C\cc$-probability space with $\tau$ a faithful state on $\CA$, and let
$\{x_1, \ldots, x_{r+s}\}$ be a semicircular system in $(\CA, \tau)$. Then there is a $P$-null set $N\subseteq\Omega$ such that for any $\omega\in \Omega\setminus N$ and for any polynomial $p$ in $r+s$ non-commuting variables:
\[
\lim_{n\rightarrow \infty}\|p(X_1^{(n)}(\omega), \ldots, X_{r+s}^{(n)}(\omega))\| = \|p(x_1, \ldots, x_{r+s})\|.
\]
\end{thm}

Definitions of the random matrix ensembles $\GSE(n, \frac 1n)$ and
$\GSES (n, \frac 1n)$ were given in the Introduction.

\begin{remark}\label{symp5}

With $(\CA, \tau)$ and $\{x_1, \ldots, x_{r+s}\}$ as in Theorem~\ref{symplectic}, let $(\CC, \phi)$ be another $C^*$-probability space with $\phi$ a faithful state on $\CC$, and let $\{z_1, \ldots, z_{4(r+s)}\}$ be a semicircular system in $(\CC, \phi)$. Then, by [VDN, Proposition~5.1.3],
\[
\bigcup_{j=1}^{r+s}\Bigg\{ \frac{z_{4j-3}+z_{4j-2}}{\sqrt 2}\;,\;  \frac{z_{4j-3}-z_{4j-2}}{\sqrt 2}\;,\; z_{4j-1},\; z_{4j}\Bigg\}
\]
is also a semicircular system in $(\CC, \phi)$. Hence, 
\[
\CM_1 = \bigcup_{j=1}^{r+s}\Bigg\{\frac{z_{4j-3}+z_{4j-2}}{\sqrt 2}\;,\;  \frac{z_{4j-3}-z_{4j-2}}{\sqrt 2}\Bigg\}
\]
is a semicircular system, 
\[
\CM_2 = \bigcup_{j=1}^{r+s}\Bigg\{\frac{z_{4j}-\i z_{4j-1}}{\sqrt 2}\Bigg\}
\]
is a circular system, and the sets $\CM_1$ and $\CM_2$ are $\ast$-free. 

For $1\leq j\leq r+s$ define $Z_j\in M_2(\CC)$ by
\[
Z_j = \frac{1}{\sqrt 2}\begin{pmatrix}
 \frac{1}{\sqrt 2}(z_{4j-3}-z_{4j-2}) & \frac{1}{\sqrt 2}(z_{4j}+\i z_{4j-1})\\
\frac{1}{\sqrt 2}(z_{4j}-\i z_{4j-1}) &  \frac{1}{\sqrt 2}(z_{4j-3}+z_{4j-2})
\end{pmatrix}.
\]
According to [VDN, Proposition~5.1.6], $\{Z_1, \ldots, Z_{r+s}\}$ is then a
semicircular system in $(M_2(\CC), \tr_2\otimes \phi)$. Furthermore,
faithfulness of $(\tr_2\otimes\phi)$ on $M_2(\CC)$ implies that there is a
state-preserving unital $\ast$-isomorphism $\Phi: C^*(\unit_\CA, x_1,
\ldots, x_{r+s})\rightarrow  C^*(\unit_2\otimes\unit_\CA, Z_1, \ldots,
Z_{r+s})$ such that $\Phi(x_j)=Z_j$, $1\leq j\leq r+s$ (cf. [V1,
Remark~1.8]). Hence, Theorem~\ref{symplectic} is proved, as soon as we have
shown that there is a $P$-null set $N\subset\Omega$ such that for any $p\in
\C\<X_1, \ldots, X_{r+s}\>$ and for every $\omega\in\Omega\setminus N$: 
\[
\lim_{n\rightarrow \infty}\|p(X_1^{(n)}(\omega), \ldots, X_{r+s}^{(n)}(\omega))\| = \|p(Z_1, \ldots, Z_{r+s})\|.
\]
\end{remark}

\vspace{.2cm}

\begin{prop}\label{symplectic2}
Let $t$, $u\in \N_0$ with $t+u\geq 1$, and for each $n\in \N$, let  $X_1^{(n)}, \ldots, X_{t+u}^{(n)}$ be stochastically independent random matrices with  $X_1^{(n)}, \ldots, X_t^{(n)}\in \GOE(n,\frac 1n)$ and $ X_{t+1}^{(n)}, \ldots, X_{t+u}^{(n)}\in \GOES(n,\frac 1n)$. 
Furthermore, let $\{x_1, \ldots, x_{t+u}\}$ be a semicircular system in a $C^*$-probability space $(\CC, \phi)$ with $\phi$ a faithful state on $\CC$, and let $\CB$ be a unital exact $C\cc$-algebra. 
Then there is a $P$-null set $N\subset \Omega$ such that for any polynomial $p$ in $t+u$ non-commuting variables with coefficients in $\CB$, and for every $\omega\in \Omega\setminus N$:
\[
\lim_{n\rightarrow \infty}\|p(X_1^{(n)}(\omega), \ldots,
X_{t+u}^{(n)}(\omega))\|_{M_n(\CB)} = \|p(x_1, \ldots, x_{t+u})\|_{\CB
  \otimes_{{\rm min}} C\cc(\CC_0)},
\]
where $\CC_0 = C\cc(\unit_\CC, x_1, \ldots ,  x_{t+u})$.
\end{prop}

\proof
This follows from Theorem~A in the same way as [HT2, Theorem~9.1] follows from [HT2, Theorem~7.1]. $\endproof$

Note that Theorem~\ref{symplectic2} applies in the case $\CB = M_2(\C)$.

\vspace{.2cm}

\textit{Proof of Theorem~\ref{symplectic}} For each $n\in\N$ we may choose
stochastically independent random matrices $Z_1^{(n)}, \ldots,
Z_{4(r+s)}^{(n)}$ such that for $1\leq j\leq r$ we have that $Z_{4j-3}^{(n)} \in \GOE(n, \frac{1}{4n})$, $Z_{4j-2}^{(n)}, Z_{4j-1}^{(n)}, Z_{4j}^{(n)}\in \GOES(n, \frac{1}{4n})$, and
\begin{eqnarray*}
X_j^{(n)} &=& \unitH\otimes Z_{4j-3}^{(n)} + \j\otimes (\i Z_{4j-2}^{(n)}) + \k\otimes (\i Z_{4j-1}^{(n)}) + \l\otimes (\i Z_{4j}^{(n)})\\
& = & \unitH\otimes Z_{4j-3}^{(n)} + 
\begin{pmatrix}
-1 & 0\\
0 & 1
\end{pmatrix} \otimes Z_{4j-2}^{(n)} +
\begin{pmatrix}
0 & \i \\
-\i & 0
\end{pmatrix} \otimes Z_{4j-1}^{(n)} +
\begin{pmatrix}
0 & -1\\
-1 & 0
\end{pmatrix} \otimes Z_{4j}^{(n)}, 
\end{eqnarray*}
and for $r+1\leq j\leq r+s$ we have that $Z_{4j-3}^{(n)} \in \GOES(n, \frac{1}{4n})$, $Z_{4j-2}^{(n)}, Z_{4j-1}^{(n)}, Z_{4j}^{(n)}\in \GOE(n, \frac{1}{4n})$, and
\begin{eqnarray*}
X_j^{(n)} &=& \unitH\otimes Z_{4j-3}^{(n)} + \j\otimes (\i Z_{4j-2}^{(n)}) + \k\otimes (\i Z_{4j-1}^{(n)}) + \l\otimes (\i Z_{4j}^{(n)})\\
& = & \unitH\otimes Z_{4j-3}^{(n)} + 
\begin{pmatrix}
-1 & 0\\
0 & 1
\end{pmatrix} \otimes Z_{4j-2}^{(n)} +
\begin{pmatrix}
0 & \i \\
-\i & 0
\end{pmatrix} \otimes Z_{4j-1}^{(n)} +
\begin{pmatrix}
0 & -1\\
-1 & 0
\end{pmatrix} \otimes Z_{4j}^{(n)}.
\end{eqnarray*}

Let $p\in\C\<X_1, \ldots, X_{r+s}\>$. Then there is a polynomial $q$ in $4(r+s)$ non-commuting variables and with coefficients in $M_2(\C)$ such that
\begin{equation}\label{symp3}
p(X_1^{(n)}, \ldots, X_{r+s}^{(n)})=q(2Z_1^{(n)}, \ldots, 2Z_{4(r+s)}^{(n)}).
\end{equation}
With $(\CC, \phi)$ as in Theorem~\ref{symplectic2}, let $\{z_1, \ldots, z_{4(r+s)}\}$ be a semicircular system in $(\CC, \phi)$. Set $\CB = M_2(\C)$, and choose $N\subset \Omega$ as in  Theorem~\ref{symplectic2}, such that for any $\omega\in \Omega\setminus N$,
\begin{equation}\label{symp4}
\lim_{n\rightarrow \infty} \|q(2Z_1^{(n)}(\omega), \ldots, 2Z_{4(r+s)}^{(n)}(\omega)))\|= \|q(z_1, \ldots, z_{4(r+s)})\|. 
\end{equation}
Clearly, $q$ may be chosen in such a way that
\begin{equation*}
\begin{split}
q(z_1&, \ldots, z_{4(r+s)})\\ 
&=p\Bigg[\Bigg(
\begin{pmatrix}
\frac 12 & 0\\
0 & \frac 12
\end{pmatrix} \otimes z_{4j-3}+
\begin{pmatrix}
-\frac 12 & 0\\
0 & \frac 12
\end{pmatrix} \otimes z_{4j-2}+
\begin{pmatrix}
0 & \frac \i2\\
-\frac \i2 & 0
\end{pmatrix} \otimes z_{4j-1}+
\begin{pmatrix}
0 & -\frac 12\\
-\frac 12 & 0
\end{pmatrix} \otimes z_{4j}\Bigg)_{1\leq j \leq r+s}\Bigg]\\ 
&=p\Bigg[\Bigg(\frac{1}{\sqrt 2}\begin{pmatrix}
 \frac{1}{\sqrt 2}(z_{4j-3}-z_{4j-2}) & \frac{1}{\sqrt 2}(z_{4j}+\i z_{4j-1})\\
\frac{1}{\sqrt 2}(z_{4j}-\i z_{4j-1}) &  \frac{1}{\sqrt 2}(z_{4j-3}+z_{4j-2})
\end{pmatrix} \Bigg)_{1\leq j \leq r+s}\Bigg].
\end{split}
\end{equation*}
Combining this identity with (\ref{symp3}) and (\ref{symp4}) we find that for every $\omega\in\Omega\setminus N$,
\[
\lim_{n\rightarrow \infty} \|p(X_1^{(n)}(\omega), \ldots, X_{r+s}^{(n)}(\omega))\| = \Bigg\|    
p\Bigg[\Bigg(\frac{1}{\sqrt 2}\begin{pmatrix}
 \frac{1}{\sqrt 2}(z_{4j-3}-z_{4j-2}) & \frac{1}{\sqrt 2}(z_{4j}+\i z_{4j-1})\\
\frac{1}{\sqrt 2}(z_{4j}-\i z_{4j-1}) &  \frac{1}{\sqrt 2}(z_{4j-3}+z_{4j-2})
\end{pmatrix} \Bigg)_{1\leq j \leq r+s}\Bigg]
\Bigg\|
\]
and, according to Remark~\ref{symp5}, this completes the proof. $\endproof$

\section{Identifying $\Lambda$ in special cases.}

The aim of this section is to identify the distribution $\Lambda$ occurring
in Lemma~\ref{distribution} and Theorem~\ref{distthm}. We shall concentrate
on the cases $S_n \in \GOE(n,\frac 1n)$,  $S_n\in \GOES(n, \frac 1n)$,
$S_n\in \GSE(n,\frac 1n)$ and $S_n\in \GSES(n,\frac 1n)$. 

We let $x$ be a semicircular element in a $C^*$-probability space $(\CA, \tau)$ with $\tau$ a faithful state on $\CA$.

\vspace{.2cm}

\begin{thm}
Let $\phi\in\Ccinf$. Then 
\begin{itemize}
\item[(i)] for $n\in\N$ and $X_n\in\GOE(n,\frac 1n)$,
\begin{equation}
\E\{\tr_n(\phi(X_n))\}=\tau(\phi(x))+\frac{1}{2n}\Bigg(\frac{\phi(-2)+\phi(2)}{2}-\frac 1\pi \int_{-2}^{2}\frac{\phi(x)}{\sqrt{x^2-4}}\ \d x \Bigg) + O\Big(\frac{1}{n^2}\Big),
\end{equation}

\item[(ii)] for $n\in\N$ and $X_n\in\GOES(n,\frac 1n)$,
\begin{equation}
\E\{\tr_n(\phi(X_n))\}=\tau(\phi(x))+\frac{1}{2n}\Bigg(\phi(0)-\frac 1\pi
\int_{-2}^{2}\frac{\phi(x)}{\sqrt{x^2-4}}\ \d x \Bigg) +
O\Big(\frac{1}{n^2}\Big),
\end{equation}
\item[(iii)] for $n\in \N$ and $X_n\in \GSE(n,\frac 1n)$,
  \begin{equation}\label{dist7}
    \E\{(\tr_2\otimes\tr_n)\phi(X_n)\} = \tau(\phi(x)) +
    \frac{1}{4n}\Bigg(\frac 1\pi \int_{-2}^{2}\frac{\phi(x)}{\sqrt{x^2-4}}\
    \d x -\frac{\phi(-2)+\phi(2)}{2}\Bigg) +  O\Big(\frac{1}{n^2}\Big),
\end{equation}
and
\item[(iv)] for $n\in \N$ and $X_n\in \GSES(n,\frac 1n)$,
  \begin{equation}\label{GSES}
    \E\{(\tr_2\otimes\tr_n)\phi(X_n)\} = \tau(\phi(x)) + \frac{1}{4n}\Bigg(\frac 1\pi \int_{-2}^{2}\frac{\phi(x)}{\sqrt{x^2-4}}\
    \d x -\phi(0)\Bigg) +  O\Big(\frac{1}{n^2}\Big) .
\end{equation}

\end{itemize}
\end{thm}

\proof
(i) Let $X_n\in\GOE(n, \frac 1n)$, and for $\lambda\in \C\setminus\R$ put
\[
G_n(\lambda) = \E\{\tr_n[(\lambda\unit_n-X_n)^{-1}]\}.
\]
In this first case corresponding to $m=1$, $r=1$, $a_0=0$, and $a_1 =\frac{1}{\sqrt2}$ we have:
\begin{eqnarray*}
R(\lambda) &=& \tau[(\lambda\unit_\CA-x)^{-2}]=-g'(\lambda),\\
L(\lambda) &=& \tau[(\lambda\unit_\CA-x)^{-1}R(\lambda)g(\lambda)^{-1}(\lambda\unit_\CA-x)^{-1}]\\
& = & -R(\lambda)g(\lambda)^{-1}g'(\lambda)\\
& = & [g'(\lambda)]^2 g(\lambda)^{-1},
\end{eqnarray*}
for $\lambda\in\C\setminus[-2,2]$.

Then, according to Lemma~\ref{distribution}, 
\begin{equation}\label{identify1}
\Lambda(\phi) = \lim_{y\rightarrow 0^+}\frac{\i}{2\pi}\int_\R \phi(x)[L(x+\i y)-L(x-\i y)]\d x, \;\;\;\;\;\;(\phi\in \Ccinf)
\end{equation}
defines a distribution $\Lambda$ on $\R$ with $\supp(\Lambda)\subseteq [-2,2]$.

By [VDN, Example~3.4.2]
\[
g(\lambda) = \frac{\lambda-\lambda\sqrt{1-\frac{4}{\lambda^2}}}{2},
\]
where $\sqrt{1-\frac{4}{\lambda^2}}$ means the principal value of
$\sqrt{1-\frac{4}{\lambda^2}}$ , i.e. $\sqrt{1-\frac{4}{\lambda^2}}>0$ if $\lambda\in\R$, $\lambda >2$. Hence,
\[
g'(\lambda) = \frac 12\Bigg(1 -\frac{\lambda}{\sqrt{\lambda^2-4}}\Bigg),
\]
and 
\begin{eqnarray*}
L(\lambda) &=&  \frac 14\Bigg(1 -\frac{\lambda}{\sqrt{\lambda^2-4}}\Bigg)^2\frac{2}{\lambda-\sqrt{\lambda^2-4}}\\
&=&\frac 12 \Bigg(\frac{\lambda}{\lambda^2-4}-\frac{1}{\sqrt{\lambda^2-4}}\Bigg).
\end{eqnarray*}

With
\begin{equation}
\nu_1 = \frac 12\ (\delta_{-2}+\delta_2),
\end{equation}
$\frac{\lambda}{\lambda^2-4}$ is the Stieltjes transform of $\nu_1$, i.e.
\begin{equation}
G_{\nu_1}(\lambda) = \frac{\lambda}{\lambda^2-4}, \;\;\;\;\;\;(\lambda \in\C\setminus\{\pm 2\}).
\end{equation}
Moreover, one may calculate the moments of the probability measure $\nu_2$ given by
\begin{equation}
\d \nu_2(x) = \frac{1}{\pi} \frac{1}{\sqrt{x^2-4}}\cdot 1_{[-2,2]}(x)\d x
\end{equation}
and use the series expansion of $(\lambda^2-4)^{-1/2}$, $(|\lambda|>2)$, to see that $\frac{1}{\sqrt{\lambda^2-4}}$ is the Stieltjes transform of $\nu_2$, i.e.
\begin{equation}
G_{\nu_2}(\lambda) = \frac{1}{\sqrt{\lambda^2-4}}, \;\;\;\;\;\;(\lambda\in\C\setminus[-2,2]).
\end{equation}
Altogether
\begin{equation}
L(\lambda) = \frac 12\ (G_{\nu_1}(\lambda)-G_{\nu_2}(\lambda)),\;\;\;\;\;\;(\lambda\in\C\setminus[-2,2]).
\end{equation}
Hence, by the inverse Stieltjes transform and (\ref{identify1}), $\Lambda$ is the distribution on $\R$ corresponding to the signed measure $\frac 12(\nu_1-\nu_2)$. 

By Theorem~\ref{distthm}, for any $\phi\in\Ccinf$
\[
\E\{\tr_n(\phi(X_n))\}=\tau(\phi(x))+\frac 1n \Lambda(\phi) + O\Big(\frac{1}{n^2}\Big),
\]
that is
\[
\E\{\tr_n(\phi(X_n))\}=\tau(\phi(x))+\frac{1}{2n}\Bigg(\frac{\phi(-2)+\phi(2)}{2}-\frac
1\pi \int_{-2}^{2}\frac{\phi(x)}{\sqrt{x^2-4}}\ \d x \Bigg) + O\Big(\frac{1}{n^2}\Big).
\]

(ii) Let $X_n\in \GOES(n, \frac 1n)$, and for $\lambda\in \C\setminus\R$ put
\[ 
G_n(\lambda) = \E\{\tr_n[(\lambda\unit_n-X_n)^{-1}]\}.
\] 
In this case (corresponding to $m=1$, $r=1$, $a_0=0$, and $a_1=\frac{1}{\i\sqrt2}$),
\begin{eqnarray*}
R(\lambda) &=&-\tau[(\lambda\unit_\CA+x)^{-1}(\lambda\unit_\CA -x)^{-1}]\\
&=& -\frac{\tau[(\lambda\unit_\CA+x)^{-1}+(\lambda\unit_\CA -x)^{-1}]}{2\lambda}\\
&=& - \frac{\tau[(\lambda\unit_\CA -x)^{-1}]}{\lambda} \\
&=& - \frac{g(\lambda)}{\lambda},
\end{eqnarray*}
where the third equality follows from the fact that $x$ and $-x$ have the same distribution. Thus, 
\begin{eqnarray*} 
L(\lambda) & = & \tau[(\lambda\unit_\CA-x)^{-1}R(\lambda)g(\lambda)^{-1}(\lambda\unit_\CA-x)^{-1}]\\
& = & - \frac{\tau[(\lambda\unit_\CA -x)^{-2}]}{\lambda}\\
& = & \frac{g'(\lambda)}{\lambda}\\
& = & \frac 12 \Bigg(\frac 1\lambda - \frac{1}{\sqrt{\lambda^2 -4}}\Bigg)\\
& = & \frac 12\ (G_{\nu_3}(\lambda)-G_{\nu_2}(\lambda)),
\end{eqnarray*}
where $\nu_2$ is the measure defined above, and $\nu_3 = \delta_0.$

Hence, for any $\phi\in\Ccinf$ we have:
\[
\E\{\tr_n(\phi(X_n))\}=\tau(\phi(x))+\frac{1}{2n}\Bigg(\phi(0)- \frac 1\pi \int_{-2}^{2}\frac{\phi(x)}{\sqrt{x^2-4}}\ \d x\Bigg) + O\Big(\frac{1}{n^2}\Big).
\]

(iii) Consider a random matrix $X_n\in\GSE(n,\frac 1n)$. For convenience we introduce some notation:
\[
\CJ = \j, \;\;\; \CK = \k, \;\;\; \CL = \l,
\]
and we let $\H^\C$ denote the complexification of the quaternions. Then
\[
\H^\C \cong {\rm span}_\C\{\unit_2, \CJ, \CK, \CL\} = M_2(\C).
\]
It follows from the definition of $\GSE(n,\frac 1n)$ and the relations
between $\GOE(n,\frac 1n)$, $\GOES(n,\frac 1n)$ and $\GRMR(n,\frac 1n)$ mentioned in the Introduction, that there are independent random matrices $Y_1^{(n)}, Y_2^{(n)},
Y_3^{(n)}$ and $Y_4^{(n)}$ from $\GRMR(n,\frac 1n)$ such that
\[
X_n = \sum_{j=1}^4  (a_j\otimes Y_j^{(n)}+ a_j^*\otimes {Y_j^{(n)}}^*),
\]
where
\[
a_1 = \frac{1}{2\sqrt2}\ \unit_2, \;\;\; a_2 =  \frac{1}{2\sqrt2}\ \CJ, \;\;\;
a_3=\frac{1}{2\sqrt2}\ \CK, \;{\rm and}\; a_4 =  \frac{1}{2\sqrt2}\ \CL.
\]
Let $(\CB, \tau)$ be any $C^*$-probability space with $\tau$ a faithful
state on $\CB$, and let $\{y_1, y_2, y_3, y_4\}$ be a circular system in
$(\CB, \tau)$. Define $s\in \H^\C\otimes\CB$ by
\[
s= \sum_{j=1}^4  (a_j\otimes y_j+ a_j^*\otimes y_j^*).
\]
Then
\begin{eqnarray*}
s & = & \frac{1}{\sqrt2}
\begin{pmatrix} 
 \frac{1}{\sqrt2}\Big( \frac{1}{\sqrt2}(y_1+y_1\cc)+ \frac{\i}{\sqrt2}(y_2-y_2\cc)\Big) & \frac{1}{\sqrt2}\Big( \frac{1}{\sqrt2}(y_3-y_3\cc)+ \frac{\i}{\sqrt2}(y_4-y_4\cc)\Big)\\
 \frac{1}{\sqrt2}\Big( \frac{1}{\sqrt2}(y_3\cc-y_3)+ \frac{\i}{\sqrt2}(y_4-y_4\cc)\Big) &  \frac{1}{\sqrt2}\Big( \frac{1}{\sqrt2}(y_1+y_1\cc)+ \frac{\i}{\sqrt2}(y_2\cc-y_2)\Big)
\end{pmatrix}\\
& = &  \frac{1}{\sqrt2} \begin{pmatrix}
 \frac{1}{\sqrt2}(x_1+x_2) &  \frac{1}{\sqrt2}(\i x_3 +x_4)\\
 \frac{1}{\sqrt2}(-\i x_3 +x_4) & \frac{1}{\sqrt2}(x_1-x_2)
\end{pmatrix}
\end{eqnarray*}
for a semicircular system $\{x_1, \ldots, x_4\}$ in $(\CB, \tau)$. By  [VDN, Proposition~5.1.3] and [VDN, Proposition~5.1.6], this representation of $s$ reveals that $s$ is circular in $(M_2(\CB), \tr_2\otimes \tau)$. Hence, it suffices to prove that
\begin{equation}\label{dist9}
    \E\{(\tr_2\otimes\tr_n)\phi(X_n)\} = (\tr_2\otimes\tau)\phi(s) +
    \frac{1}{4n}\Bigg(\frac 1\pi \int_{-2}^{2}\frac{\phi(x)}{\sqrt{x^2-4}}\
    \d x -\frac{\phi(-2)+\phi(2)}{2}\Bigg) +  O\Big(\frac{1}{n^2}\Big).
\end{equation}

According to Theorem~\ref{distthm} there is a distribution $\Lambda$ with
$\supp(\Lambda)\subseteq \sigma(s)=[-2,2]$, such that
\[
\E\{(\tr_2\otimes\tr_n)\phi(X_n)\}=(\tr_2\otimes\tau)\phi(s)+\frac 1n \Lambda(\phi) + O\Big(\frac{1}{n^2}\Big),
\]

As in the two first cases, to identify $\Lambda$ we try to recognize $l(\lambda)=\tr_2(L(\lambda\unit_2))$ as the Stieltjes
transform of a signed measure and then apply the inverse Stieltjes transform. By definition
\begin{equation}\label{dist6}
l(\lambda) = (\tr_2\otimes\tau) [(\lambda\unit_2\otimes\unit_\CB -
s)^{-1}(R(\lambda\unit_2)G(\lambda\unit_2)^{-1}\otimes\unit_\CB)(\lambda\unit_2\otimes\unit_\CB -
s)^{-1}],
\end{equation}
where (note that in this case $a_j^*= \pm a_j$)
\[
R(\lambda\unit_2) = 2\sum_{j=1}^4\sum_{k,l=1}^2(\id_2\otimes\tau)[(a_j
e_{kl}^{(2)}\otimes\unit_\CB)(\lambda\unit_2\otimes\unit_\CB - \overline{s})^{-1}(e_{kl}^{(2)}a_j\otimes\unit_\CB)(\lambda\unit_2\otimes\unit_\CB -
s)^{-1}].
\]
We prove that
\begin{equation}\label{dist1}
  R(\lambda\unit_2) = -\frac 12 \ (\id_2\otimes\tau)[(\lambda\unit_2\otimes\unit_\CB
  - s)^{-2}] = \frac 12\ \diffl G(\lambda\unit_2).
\end{equation}

To this end define a linear map $\Psi : \H^\C\otimes\CB \rightarrow
\H^\C\otimes\CB$ by
\begin{eqnarray*}
  \Psi(\unit_2\otimes b) &=& \unit_2\otimes b, \\
  \Psi(\CJ\otimes b) &=& -\CJ\otimes b, \\
  \Psi(\CK\otimes b) &=& \CK\otimes b, \\
  \Psi(\CL\otimes b) &=& -\CL\otimes b,
\end{eqnarray*}
for $b\in \CB$. One easily checks that $\Psi$ is actually a (well-defined) unital $\ast$-isomorphism. In particular, as
\[
\lambda\unit_2\otimes\unit_\CB - \overline{s} =
\Psi(\lambda\unit_2\otimes\unit_\CB -s),
\]
we have that
\begin{equation}\label{dist2}
(\lambda\unit_2\otimes\unit_\CB - \overline{s})^{-1} = (\Psi(\lambda\unit_2\otimes\unit_\CB -s))^{-1} = \Psi((\lambda\unit_2\otimes\unit_\CB -s)^{-1}).
\end{equation}

Let $x\otimes v$ be an elementary tensor in $\H^\C\otimes\CB$. Then
\[
  \sum_{k,l=1}^2(e_{kl}^{(2)}\otimes\unit_\CB)\Psi(x\otimes
  v)(e_{kl}^{(2)}\otimes\unit_\CB) =
  \begin{cases}
    x\otimes v, &{\rm if} \;\;x= \unit_2,\\
    -x\otimes v, &{\rm if} \;\;x= \CJ,\\
    -x\otimes v, &{\rm if} \;\;x= \CK,\\
    -x\otimes v, &{\rm if} \;\;x= \CL,
  \end{cases}
\]
and it is standard to check that this implies that for $x\in\{\unit_2, \CJ, \CK, \CL\}$ (and hence for any $x\in\H^\C$) we have:
\begin{equation}\label{dist3}
\sum_{j=1}^4\sum_{k,l=1}^2 (a_j
e_{kl}^{(2)}\otimes\unit_\CB)\Psi(x\otimes
v)(e_{kl}^{(2)}a_j\otimes\unit_\CB) =
   -\frac 14 (x\otimes v).
\end{equation}

As $\lambda\unit_2\otimes\unit_\CB - \overline{s}\ $ belongs to the unital $C^*$-algebra $\H^\C\otimes\CB$, $(\lambda\unit_2\otimes\unit_\CB - \overline{s})^{-1}$ also belongs to $\H^\C\otimes \CB$. Thus, (\ref{dist2}) and  (\ref{dist3}) imply that
\[
\sum_{j=1}^4\sum_{k,l=1}^2(a_je_{kl}^{(2)}\otimes\unit_\CB)(\lambda\unit_2\otimes\unit_\CB
- \overline{s})^{-1} (e_{kl}^{(2)}a_j\otimes\unit_\CB) = -\frac 14
\ (\lambda\unit_2\otimes\unit_\CB -s)^{-1},
\]
and consequently,
\begin{equation*}
R(\lambda\unit_2) = -\frac 12 \ (\id_2\otimes\tau)[(\lambda\unit_2\otimes\unit_\CB
  - s)^{-2}] = \frac 12 \ \diffl G(\lambda\unit_2).
\end{equation*}

The next step is to prove that $G(\lambda\unit_2)\in\C\unit_2$. We have seen that
\[
G(\lambda\unit_2)^{-1} = \lambda\unit_2-\sum_{j=1}^4(a_jG(\lambda\unit_2)a_j\cc + a_j\cc G(\lambda\unit_2)a_j),
\]  
that is
\begin{equation}\label{dist4}
G(\lambda\unit_2)^{-1}  =  \lambda\unit_2 - \frac 14 \ (G(\lambda\unit_2)-\CJ G(\lambda\unit_2)\CJ - \CK G(\lambda\unit_2)\CK - \CL G(\lambda\unit_2)\CL).
\end{equation}
Now, for any $x\in M_2(\C)$ we have:
\begin{equation}\label{dist8}
x-\CJ x\CJ - \CK x\CK - \CL x\CL = 4\ \tr_2(x)\unit_2 \in \C \unit_2,
\end{equation}
and therefore, by (\ref{dist4}), $G(\lambda\unit_2)^{-1}\in\C\unit_2$. Then $G(\lambda\unit_2)$ is a scalar too, so 
\[
G(\lambda\unit_2) = \tr_2( G(\lambda\unit_2))\ \unit_2 = g(\lambda)\ \unit_2,
\]
and
\[
G(\lambda\unit_2)^{-1} = g(\lambda)^{-1}\ \unit_2.
\]
By (\ref{dist1})
\[
R(\lambda\unit_2) = \frac 12 \ g'(\lambda)\ \unit_2,
\]
and again, since $s$ is semicircular, 
\[
g(\lambda) = \frac{\lambda-\lambda\sqrt{1-\frac{4}{\lambda^2}}}{2}.
\]

Inserting these expressions into (\ref{dist6}) we find that for $\lambda\in\C\setminus [-2,2]$ we have:
\begin{eqnarray*}
l(\lambda)& = & \frac 12 \ g'(\lambda) \ g(\lambda)^{-1}\ (\tr_2\otimes\tau)[(\lambda\unit_2\otimes\unit_\CB-s)^{-2}] \\
& = & -\frac 12\ [g'(\lambda)]^2 \ g(\lambda)^{-1}\\
& = & -\frac 12\ \frac 14 \Bigg(1- \frac{\lambda}{\sqrt{\lambda^2-4}}\Bigg)^2 \frac{2}{\lambda - \sqrt{\lambda^2 -4}}\\
& = & \frac 14 \ \Bigg(\frac{1}{\sqrt{\lambda^2 -4}}-\frac{\lambda}{\lambda^2-4}\Bigg)\\
& = & \frac 14\ (G_{\nu_2}(\lambda) - G_{\nu_1}(\lambda)),
\end{eqnarray*}
where $\nu_1$ and $\nu_2$ are the measures defined above. Then
(\ref{dist9}) follows as in the previous cases.

(iv) The proof of (iv) is similar to the proof of (ii), but one must apply
some of the techniques from the proof of (iii) too. We leave out the details. $\endproof$

\vspace{.5cm}

{\bf Acknowledgments.} I would like to thank my advisor, Uffe~Haagerup, with whom I had many enlightening discussions, and who
made some important contributions to this paper. Also, thanks to Steen~Thorbj\o rnsen who took time to answer several questions. 

{\small

\noindent Department of Mathematics and Computer Science\\
University of Southern Denmark\\
Campusvej~55, 5230~Odense~M\\
Denmark\\
{\tt schultz@imada.sdu.dk}

\end{document}